\documentclass[11pt]{article}

\pdfoutput=1
\usepackage[pdftex]{color}
\usepackage{amssymb}
\usepackage{amsthm}
\usepackage{amsmath}
\usepackage{latexsym}
\usepackage{amscd}
\usepackage{graphicx}
\usepackage[pdftex, colorlinks=true, citecolor=green]{hyperref}
\usepackage{lscape}
\usepackage{multirow}

\newcommand{\ds}{\displaystyle}

\newcommand{\ben}{\begin{equation}}     
\newcommand{\eeqn}{\end{equation}}
\newcommand{\bey}{\begin{eqnarray}}
\newcommand{\eey}{\end{eqnarray}}

\newtheorem{thm}{Theorem}[section]
\newtheorem{prop}[thm]{Proposition}
\newtheorem{lemma}[thm]{Lemma}

\newtheorem{defn}[thm]{Definition}
\newtheorem{conj}[thm]{Conjecture}

\begin{document}

\begin{flushleft}
{\Large
\textbf{Networks of coupled quadratic nodes}
}
\\
\vspace{4mm}
 Anca R\v{a}dulescu$^{*,}\footnote{Assistant Professor, Department of Mathematics, State University of New York at New Paltz; New York, USA; Phone: (845) 257-3532; Email: radulesa@newpaltz.edu}$, Simone Evans$^1$
\\
\indent $^1$ Department of Mathematics, SUNY New Paltz, NY 12561
\\
\end{flushleft}

\vspace{3mm}
\begin{abstract}
We study asymptotic dynamics in networks of coupled quadratic nodes. While single map complex quadratic iterations have been studied over the past century, considering ensembles of such functions, organized as coupled nodes in a network, generate new questions with potentially interesting applications to the life sciences. 

We investigate how traditional Fatou-Julia results may generalize in the case of networks. We discuss extensions of concepts like escape radius, Julia and Mandelbrot sets (as parameter loci in $\mathbb{C}^n$, where $n$ is the size of the network). We study topological properties of these asymptotic sets and of their two-dimensional slices in $\mathbb{C}$ (defined in previous work). We find that, while network Mandelbrot sets no longer have a hyperbolic bulb structure, some of their geometric landmarks are preserved (e.g., the cusp always survives), and other properties (such as connectedness) depend on the network structure. We investigate possible extensions of the relationship between the Mandelbrot set and the Julia set connectedness loci in the case of network dynamics. 

We discuss possible classifications of asymptotic behavior in networks based on their underlying graph structure, using the geometry of Julia or Mandelbrot sets as a classifier. Finally, we propose a method for book-keeping asymptotic dynamics simultaneously over many networks with a common graph-theoretical property. \emph{Core} Julia and Mandelbrot sets describe statistically average asymptotic behavior of orbits over an entire collection of configurations.
\end{abstract}

\section{Introduction}
\label{introd}

\subsection{Networks of complex quadratic maps}

Many natural systems are organized as self-interacting networks. Subsequently, dynamic networks have been used as a modeling framework in many fields of life sciences, with the definition of nodes and edges depending on the context. When studying brain networks, the nodes may represent neurons, and the connecting oriented edges between them are synapses with varying weights. In epidemics, the nodes may be populations, and the edges, the physical contacts that promote contagion. For a traffic map, the nodes may be towns, connected by various size roads, and for a social network, individuals are connected by friendship edges of different strengths. A unifying questions for all these different fields regards how the hardwired \emph{structure} of a network (its underlying graph) and its connectivity (edge weights) affect the system's overall \emph{function}.

When translating connectivity patterns to network dynamics, the main difficulty reportedly arises from the graph complexity compounded with the nodes' dynamic richness. In order to better understand this dependence, we started to investigate it in simple theoretical models, where one may more easily identify and pair specific structural patterns to their effects on dynamics. Our choice is further motivated by the fact that, historically, discrete iterations have provided good simplified representations for many natural processes such as learning in brain circuits.
 
We are focused in particular on understanding how architecture affects asymptotic dynamics in networks of complex quadratic maps. That is because, historically, the classical theory for single function iterations has been most developed for the complex quadratic family: $f_c \colon \mathbb{C} \to \mathbb{C}, \; f_c(z) = z^2+c$, for $c \in \mathbb{C}$. Work on this family spans more than a century, from the original results of  Fatou and Julia, describing in the early 1900s the behavior of orbits in the dynamic complex plane (reflected by the structure of the Julia set)~\cite{julia1918memoire,fatou1919equations}, to bifurcation phenomena in the parameter plane (reflected in the work of Mandelbrot and others, in the 1970s)~\cite{mandelbrot1980fractal,brooks1981dynamics}, to recent connections between the two concepts~\cite{branner1992iteration,bonifant2010cubic,qiu2009proof}. Therefore, we adopted the simplified framework of networked logistic maps as an ideal starting point for approaching basic dynamic questions in the context of networks. In this framework, each network node receives weighted inputs from the adjacent nodes, and integrates these inputs in discrete time as a complex quadratic map. Then the system takes the form of an iteration in $\mathbb{C}^n$:

\begin{eqnarray}
z_j(t) \longrightarrow z_j (t+1) &=& f_j\left(\sum_{k=1}^{n}{g_{jk} A_{jk} z_k} \right) \nonumber 
\label{mothermap}
\end{eqnarray}

\noindent where $n$ is the size of the network, $\ds A=(A_{jk})_{j,k =1}^n$ is the binary adjacency matrix of the oriented underlying graph, and $g_{jk}$ are the weights along the adjacency edges. In isolation, each node $z_j(t) \to z_j(t+1)$, $1 \leq j \leq n$, iterates as a quadratic function $f_j(z) = z^2 + c_j$. When coupled as a network with adjacency $A$, each node will act as a quadratic modulation on the sum of the inputs received along the incoming edges (as specified by the values of $A_{jk}$, for $1 \leq k \leq n$).

In our proposed work, we will use properties of multi-dimensional orbits in $\mathbb{C}^n$, in particular their asymptotic behavior (via the topological and fractal structure of Julia and Mandelbrot multi-sets) -- to classify dynamic behavior for different network architectures. By imposing additional structural conditions on edge density or distribution, we will investigate whether it is possible to predict the geometry of Julia and Mandelbrot sets from specific information on the network hardwiring. We aim to tease apart the instances in which small perturbations in the position or strength of one single connection may lead to dramatic topological changes in the asymptotic sets, from the instances in which these sets are robust to much more significant changes. 

In our previous work~\cite{radulescu2017real}, when suggesting possible ties of our results with broader applications to the life sciences, we interpreted iterated orbits as describing the temporal evolution of an evolving system (e.g., learning neural network). An escaping initial condition (whether in the complex plane $\mathbb{C}$, for a single iterated map, or in $\mathbb{C}^n$, for an iterated network) may be seen as an eventually unsustainable feature of the system, while a prisoner may represent a trivial, or inefficient feature. The Julia set is formed of all the boundary points between prisoners and escapees, hence we suggested that it can be regarded as the ``critical locus'' of states with a complex temporal evolution, characteristic to living systems operating within an optimal range.

In  a previous paper we defined the network Mandelbrot set, for simplicity, as the node parameter range for which the critical point (i.e., all nodes equal zero) is bounded (i.e., functionally sustainable) under iterations of the network. In the traditional case of a single iterated quadratic map, this is equivalent to defining the parameter locus for which the Julia set is connected. Indeed, the Fatou-Julia Theorem delivers in this case a well-known duality: a bounded critical orbit implies a connected Julia set, and an escaping critical orbit implies a totally disconnected Julia set. We don't expect this equivalence to remain true when iterating networks. For networks, we have already noticed that the situation is a lot more complicated: the Julia set may not necessarily be connected or totally disconnected, and may have a finite number of connected components. What we conjectured, in a slightly different form, is that a connected ``uni-Julia set'' implies a bounded critical orbit, but not conversely.

One may further interpret that, in the case of a network with connected Julia set, all sustainable initial conditions (i.e., prisoners, of initial points leading to bounded orbits) can be reached by perturbations from rest (i.e. from the critical point, with all nodes set at zero), without having to leave the prisoner set. Totally disconnected Julia sets represent a scattered, measure zero locus of sustainable initial states. We further conjectured that one would always have to traverse an intermediate asymptotic region characterized by disconnected Julia sets when transitioning from the parameters locus for connected Julia sets to the parameter locus for totally disconnected Julia sets.

\subsection{Prior results in small networks}
\label{prior}

In order to establish a conceptual framework, in previous work we considered simple, low-dimensional networks, which are both analytically tractable and allow easy visualization and interpretation of the results, suggesting a baseline for extensions to higher dimensional, more complex networks. 

We considered in particular three dimensional networks with various coupling geometries between their complex nodes $z_1$, $z_2$, $z_3$. For fixed  logistic parameters $c_1$=$c_2$=$c_3$=$c$, we described the dependence of the Julia and Mandelbrot sets and of their one-dimensional slices on the graph wiring and of the strengths of the connections between nodes.

Our prior work~\cite{radulescu2017real} suggests that even basic results from the case of a single iterated quadratic maps may have to be rediscovered in the context of networks (one yet needs to prove, for example, even the existence of an escape radius). In our study of dynamics in small quadratic networks, we redefined extensions of some of the traditional concepts: multi-orbits, Julia and Mandelbrot sets, as well as their one-dimensional complex slices, which we called uni-Julia and equi-Mandelbrot sets.

\begin{defn}
For a fixed parameter $(c_1,...,c_n) \in \mathbb{C}^n$, we call the \textbf{prisoner set} of the network, the locus of $(z_1,...,z_n) \in \mathbb{C}^n$ which produce a bounded multi-orbit in $\mathbb{C}^n$.  We call the \textbf{uni-prisoner set}, the locus of $z \in \mathbb{C}$ so that $(z,...z) \in \mathbb{C}^n$ produces a bounded multi-orbit. The \textbf{multi-Julia set (or the multi-J set)} of the network is defined as the boundary in $\mathbb{C}^n$ of the multi-prisoner set. Similarly, one defines the \textbf{uni-Julia set (or uni-J set)} of the network as the boundary in $\mathbb{C}$ of the uni-prisoner set for that network.
\end{defn}

\begin{defn}
We define the \textbf{multi-Mandelbrot set (or the multi-M set)} of the network the parameter locus of $(c_1,...,c_n) \in \mathbb{C}^n$ for which the multi-orbit of the critical point $(0,...,0)$ is bounded in $\mathbb{C}^n$. We call the \textbf{equi-Mandelbrot set (or the equi-M set)} of the network, the locus of $c \in \mathbb{C}$ for which the critical multi-orbit is bounded for \textbf{equi-parameter} $(c_1,c_2,...c_n)=(c,c,...c) \in \mathbb{C}^n$.  We call the \textbf{$k$th node equi-M set} the  locus $c \in \mathbb{C}$ such that the component of the multi-orbit of $(0,...,0)$ corresponding to the $k$th node remains bounded in $\mathbb{C}$.
\end{defn}

\noindent With these definitions, we pointed out new, network phenomena, and proposed new versions of the traditional theorems for the case of networked nodes. We showed that even in networks where all nodes are identical maps, their behavior may not be ``synchronized,'' in the sense that different nodes may have different asymptotic behavior (reflected in differences between node-wise Mandelbrot and Julia sets). Node coupling seems to enhance this ``de-synchronization'' between two or more nodes, and additional networking may generally lead to smaller network Mandelbrot and Julia sets. Unlike for the traditional, single map iterations, the definition requirement for the M-set that the origin has a bounded multi-orbit is no longer equivalent with that of the J-set being connected, in either of its forms (multi-J or uni-J set). In our previous work, however, we have conjectured a weaker version of the Fatou-Julia theorem in this case, which remains to be verified analytically. We also analyzed and interpreted the distinct effects of varying excitatory versus inhibitory strength, and those of introducing feedback into the network.

We finally pointed out that complex natural networks are typically a lot larger than the three and four node networks we had studied. At the same time, however, natural networks (such as brain circuits, for example) tend to be highly hierarchic, with the behavior of each one node at a certain complexity level integrating the behavior of a collection of lower-level nodes. Hence, at each complexity level, the size of the network to be studied may be in fact relatively small (tens or hundreds of nodes). While for small networks the effects of architecture on asymptotic dynamics can still be observed and studied by looking at each configuration individually, and for very large networks one may take the large size limit approach traditional in random graph theory, for these intermediate networks one has to build a different approach. One possible framework is statistical. Using book-keeping methods developed in our previous work~\cite{radulescu2015nonlinear}, we define probabilistic (or average) versions of the Julia and Mandelbrot sets, illustrating the likelihood that each initial state of the network remains bounded when iterated under a random network configuration with certain given properties. Using this framework, one can attempt to tease apart graph theoretical features (e.g.,  hubs, motifs) determinant of certain dynamics of the network, from those less consequential to temporal behavior.\\

\noindent Within the current paper, we will focus on studying equi-M and uni-J sets for networks with identical nodes (i.e., network dynamics for equi-parameters $c \in \mathbb{C}$).  The paper is organized as follows. In Section~\ref{extensions}, we state sufficient conditions for existence of an escape radius, and we calculate this radius in terms of the network parameters. We investigate an example family of low-dimensional networks, exploring topological properties of equi-Mandelbrot and uni-Julia sets. In Section~\ref{networks}, we introduce new methods applicable to higher dimensional networks. We investigate the robustness of the asymptotic sets under changes in the graph structure, and explore classifications. We introduce average (``core'') Julia and Mandelbrot sets, as a bookkeeping approach to simultaneously recording the properties of these sets for many configurations, in the case of higher dimensional networks. Finally, in Section~\ref{discussion}, we interpret our results and present some potential applications.


\section{Extensions of traditional results}
\label{extensions}

\subsection{Escape radius}
\label{escape}

Suppose that we have a network in which all nodes act nontrivially onto themselves (that is, each node $z_j$ has a self-loop of weight $g_{jj} \neq 0$). Moreover, suppose that this self-action is in each node larger than the sum of the outside inputs: $\lvert g_{jj} \rvert > \sum_{l \neq j} \lvert g_{jl} \rvert$, for all $j$. With this condition, we can take $\delta$ such that 
$$\frac{\lvert g_{jj} \rvert}{\sum_{l \neq j} \lvert g_{jl} \rvert} > \delta > 1$$

\noindent Then one can show that the network has escape radius $R$ that depends on the network weights.

\begin{lemma}
There exists a large enough $M$ such that, if $\lvert z_j(k) \rvert \leq M$ for all nodes $1 \leq j \leq n$ at all iterates $0 \leq k \leq K+1$, then it follows that $\ds \lvert z_j(k) \rvert \leq \frac{M}{\delta}$, for all $1 \leq j \leq n$ and all $0 \leq k \leq K$.
\label{lemmaM}
\end{lemma}

\proof{Take an $M>0$ large enough so that $\lvert z_j(k) \rvert \leq M$ for all nodes $1 \leq j \leq n$ at all iterates $0 \leq k \leq K+1$. Recall that, for all $1 \leq j \leq n$ and for all $1 \leq k \leq K$, we have:
$$\lvert z_j(k+1) \rvert = \left \lvert \left( \sum_l g_{jl}z_l(k) \right)^2 +c_j \right \rvert \geq \left \lvert \sum_l g_{jl}z_l(k) \right \rvert^2 -\lvert c_j \rvert \; \Longrightarrow$$  
$$\left \lvert \sum_l g_{jl}z_l(k) \right \rvert \leq \sqrt{\lvert z_j(k+1) \rvert + \lvert c_j \rvert} \leq \sqrt{M+\lvert c_j \rvert} \; \Longrightarrow$$
$$\lvert g_{jj} \rvert \lvert z_j(k) \rvert \leq \sqrt{M+\lvert c_j \rvert} + \sum_{l \neq j} \lvert g_{jl} \rvert \lvert z_l(k) \rvert \leq \sqrt{M + \lvert c_j \rvert} + M\sum_{l \neq j} \lvert g_{jl} \rvert \; \Longrightarrow$$
$$\lvert z_j(k) \rvert \leq \frac{\sqrt{M + \lvert c_j \rvert} + MG_j}{\lvert g_{jj} \rvert}$$

\noindent where $G_j =\sum_{l \neq j} \lvert g_{jl} \rvert$ is the sum of the external input weight onto each node $z_j$. We ask for a sufficient condition for $M$ that would insure that the right side
$$\frac{\sqrt{M + \lvert c_j  \rvert}+MG_j}{\lvert g_{jj} \rvert} \leq \frac{M}{\delta} \; \Longrightarrow \; 2\sqrt{M+ \lvert c_j \rvert} \leq M A_j$$

\noindent where $A_j = \lvert g_{jj} \rvert - \delta \sum_{l \neq j} \lvert g_{jl} \rvert >0$ (as per our original assumption on node-wise input strengths). Then we can square both sides and get: $\ds M^2 A_j^2 - \delta^2 M -\delta^2 \lvert c_j \rvert \geq 0$, which can be easily accomplished if $M$ is taken to be larger than the higher quadratic root: 
$$M > \frac{\delta^2 + \sqrt{\delta^4 + \delta^2 \lvert c_j \rvert^2 A_j^2}}{2A_j^2}$$

}

\noindent {\bf Remark.} Under the existing assumptions on the network weights and the notation in Lemma~\ref{lemmaM}, the network dynamics has escape radius $M/\delta$. More precisely: if $z_j(k) > M/\delta$ for some node $j$ at the iteration step $k$, then $z_l(k+1) > M$ for some node $l$ at next  iteration step $k+1$. Hence, once it left the disc of radius $M/\delta$, the network will go to $\infty$ in the maximum norm. In conclusion, we have the following theorem:

\begin{thm}
If a network satisfies the condition $\lvert g_{jj} \rvert > \sum_{l \neq j} \lvert g_{jl} \rvert$ for all $j$, then it has escape radius that depends on the network weights.
\end{thm}

\subsection{Main cardioid and periodic bulbs}

Possibly the most striking geometric features of the traditional Mandelbrot set are its periodic Fatou components. Indeed, one may consider the set ${\cal M}'$ of all parameters $c$ for which the map $f_c$ has an attracting periodic orbit. It has been established that ${\cal M}'$ is a subset of the interior of ${\cal M}$. For example, the $c$-locus for which the map has an attracting fixed point represents the interior of the main cardioid of ${\cal M}$, and the locus for which the map has an attracting period two orbit is the interior of the disc of radius $1/4$ centered at $(-1,0)$. Whether ${\cal M}'$ is in fact identical to the interior of ${\cal M}$, or ${\cal M}$ contains other (``ghost'', non-hyperbolic) interior points -- is still an open question, know as the Density of Hyperbolicity conjecture. While the conjecture was solved for real polynomials over twenty years ago~\cite{lyubich1997dynamics,graczyk1997generic}, it still represents one of the most important open problems in complex dynamics.

In the traditional case of single iterated maps, the hyperbolic components (bulbs) of the Mandelbrot set are identified with the parameter subsets for which the map has an attracting orbit of period $k$. For example, the locus in $\mathbb{C}$ for which the map has an attracting fixed point is the interior of the main cardioid, defined as $\ds c = \frac{e^{i\theta}}{2} - \frac{e^{2i\theta}}{4}$, with $0 \leq \theta \leq 2\pi$.

One can similarly compute the hyperbolic components for a network of quadratic complex nodes. To fix our ideas, we calculate the main hyperbolic component (representing the locus of $c \in \mathbb{C}$ for which the network has an attracting fixed point) for a very simple network of three nodes (which we had considered in previous work). We will illustrate how the boundary of this region differs from the main cardioid from the traditional case, and compare it with the numerical illustrations of the corresponding equi-M sets.

Consider the following ``simple dual'' network with two input and one output nodes:
\begin{eqnarray*}
z_1 &\to& z_1^2 +c\\
z_2 &\to& (az_1+z_2)^2+c\\
z_3 &\to& (z_1+z_2)^2+c
\end{eqnarray*}

\noindent where $a$ is the level of cross-talk between the input nodes. We ask that $(z_1,z_2,z_3)$ be a fixed point: $z_1^2+c =z_1$ and $(az_1+z_2)^2+c = z_2$ (with this, the third component is fixed automatically, since it is independent of $z_3$). We additionally require that the fixed point be attracting. The Jacobian matrix of this network 
\begin{equation}
J(z_1,z_2,z_3) =\left(  \begin{array}{lll} 2z_1 & 0 & 0\\ 2a(az_1+z_2) & 2(az_1+z_2) & 0\\ 2(z_1+z_2) & 2(z_1+z_2) & 0 \end{array} \right)
\end{equation}

\vspace{2mm}
\noindent has eigenvalues $\lambda_1=0$, $\lambda_2 = 2z_1$ and $\lambda_3 = 2(az_1+z_2)$. To find the boundary of the hyperbolic component, we require that $\lvert \lambda_2 \rvert = \lvert \lambda_3 \rvert = 1$ at one of the fixed points, separating the region where this fixed point is stable (attracting) from the region where it has unstable (saddle) behavior. For simplicity, call $\varphi  = az_1+z_2$, and notice that the eigenvalue condition implies $2 z_1 =e^{i\theta}$, with $ 0 \leq \theta \leq 2\pi$, and $2 \varphi = e^{i\tau}$, with $0 \leq \tau \leq 2\pi$.

\begin{figure}[h!]
\begin{center}
\includegraphics[width=0.7\textwidth]{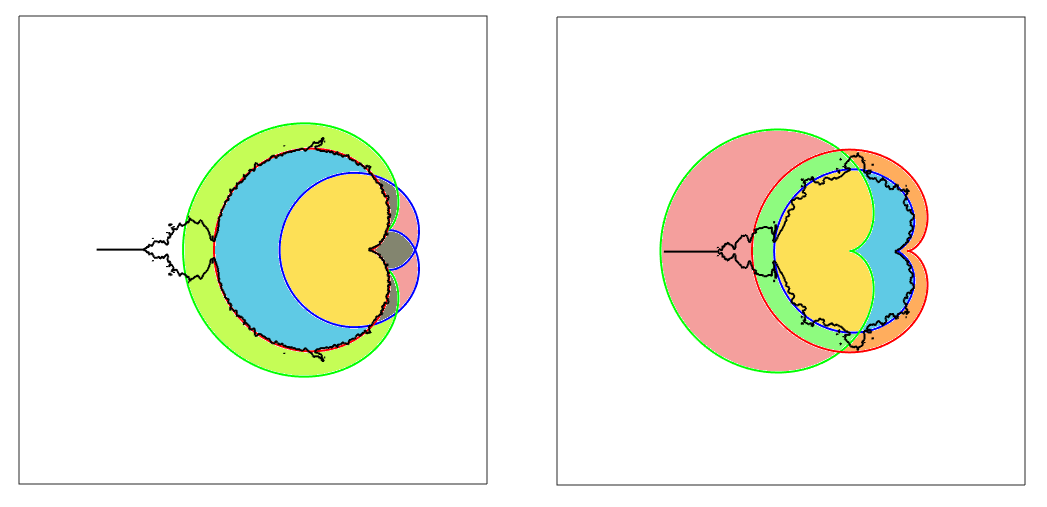}
\end{center}
\caption{\emph{\small {\bf Main hyperbolic component of the M-set, for a simple dual network with different cross-talk values.} The panels show the M-set computed numerically and the curves obtained analytically for {\bf A.} $a=-1/3$; {\bf B.} $a=1/3$. The red curve represents the traditional Mandelbrot cardioid, and the green and blue curves represent the additional restrictions for $c$, as described in the text. The color represent different stability behaviors of the critical components.}}
\label{cardio_curves}
\end{figure}

The first condition implies that $c = z_1 - z_1^2 = \frac{e^{i\theta}}{2} - \frac{e^{2i\theta}}{4}$, which is precisely the main cardioid from the traditional case. The second condition will add another restriction, which will depend on parameter $a$. It follows immediately, however, that the network hyperbolic component will always be a subset of the interior of the main cardioid from the traditional case of single map iterations.

Notice now that the first fixed point equation multiplied by $a$ and added to the second delivers: $az_1^2=\varphi-\varphi^2-(a+1)c$, while the second gives us: $\varphi^2 + c = \varphi -az_1$, hence $az_1 = \varphi-\varphi^2-c$. In conclusion:
$$ a^2z_1^2 = a\varphi - a\varphi^2 - a(a+1)c = (\varphi-\varphi^2-c)^2$$

\noindent Calling $\xi = \varphi - \varphi^2$, we obtain the quadratic equation in $c$:
$$c^2 + (a^2+a-2\xi)c + \xi^2-a\xi=0$$

\noindent which gives the solution curves:
\begin{eqnarray*}
c &=& \frac{2\xi-a-a^2 \pm \sqrt{a^2(a+1)^2-4a^2\xi}}{2} \\ \nonumber
&=&  \frac{2(\varphi-\varphi^2)-a-a^2 \pm \sqrt{a^2(a+1)^2-4a^2(\varphi-\varphi^2)}}{2} \\ \nonumber
\end{eqnarray*}

\noindent where $\varphi = e^{i\tau}/2$. We represented these curves and the regions between them in Figure~\ref{cardio_curves}.\\

\noindent Notice that  having an attracting fixed point for the network no longer implies that the origin will be in the attraction basin of this fixed point, hence the critical orbit can still escape (as shown in Figure~\ref{cardio_curves}). Hence even in networks as simple as this family of examples, structuring the interior of M-set as a union of hyperbolic bulbs fails. While some of the bulb geometry is preserved (e.g., the cusp seems robust under network transformations), some of the landmarks lose their dynamic context (e.g. the origin $c=0$, while still in the network Mandelbrot set, can no longer be regarded as the center of a main cardioid).

\begin{figure}[h!]
\begin{center}
\includegraphics[width=0.9\textwidth]{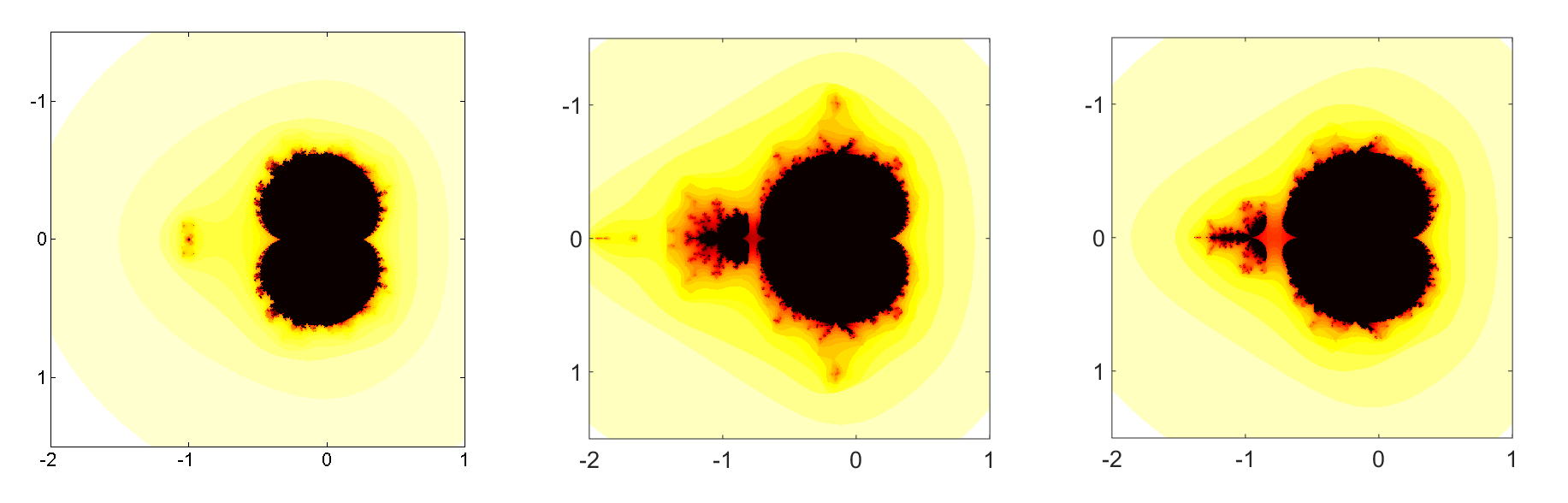}
\end{center}
\caption{\emph{\small {\bf Examples of disconnected equi-M sets for two networks with $N=3$ nodes.} All three networks belong to the family: $z_1 \to z_1^2+c$, $z_2 \to (az_1+z_2)^2+c$, $z_3 \to (z_1+z_2+bz_3)^2+c$. {\bf A.} Connectivity weights  $a = -1$, $b=-1$; {\bf B.} Connectivity weights  $a = -1/3$, $b=-1/3$; {\bf C.} with  connectivity weights  $a = -2/3$, $b=-1/3$. The colors represent the escape rate of the critical orbit out of the disc of radius $R_e=20$, so that the critical orbit is bounded in the central black region, and escapes faster with increasingly  lighter colors.}}
\label{connect}
\end{figure}

The properties of higher period bulbs get perturbed even more dramatically. We can track, for example, what becomes of the period two bulb/disc (originally centered at the $c=-1$) in a family of simple tree-dimensional networks. While part of this behavior is conserved in some networks, it completely collapses in others, depending on the configuration and connectivity parameters. To illustrate, we first show in Figure~\ref{connect} a comparison between the network Mandelbrot sets for our model network, for three different parameter pairs: $(a,b)=(-1,-1)$; $(a,b)=(-1/3,-1/3)$ and $(a,b)=(-2/3,-1/3)$. In all cases, the M-set is disconnected (see Section~\ref{Mand_connect} for detail); in the first two cases, $c=-1$ is still part of the set; in the third, it is not.

\subsection{Connectedness of the Mandelbrot set}
\label{Mand_connect}

Establishing connectedness of the traditional Mandelbrot has been historically challenging, with an original conjecture (based on numerical and visual consideration) stating the exact opposite. Connectedness of the set was finally determined by Douady and Hubbard~\cite{douady1984exploring}, with a proof based on the construction of a conformal isomorphism between the complement of the Mandelbrot set and the complement of the closed unit disk. 

It has been hypothesized that the Mandelbrot set is locally connected (the MLC conjecture). While local connectivity has been established at many special points in the Mandelbrot set (for example, Yoccoz proved that this is the case at all finitely renormalizable parameters~\cite{hubbard1992local}), the general conjecture remains open. Establishing local connectedness of the Mandelbrot is extremely desirable, since it implies Density of Hyperbolicity~\cite{douady1984exploring}.\\

\noindent It is not entirely surprising that most of these results no longer apply in this form for networked complex maps. For example, connectedness fails in general for network equi-Mandelbrot sets. To fix out ideas, we illustrate and prove disconnectedness for an example network in a three-dimensional family considered previously~\cite{radulescu2017real} (see Figure~\ref{connect}). This family (which we called the ``self-drive model'') is interesting and easy to analyze, since it is a feed-forward network (each node depends only on the ones  with smaller indices): $z_1 \to z_1^2+c$, $z_2 \to (az_1+z_2)^2+c$, $z_3 \to (z_1+z_2+bz_3)^2+c$.

\begin{figure}[h!]
\begin{center}
\includegraphics[width=0.4\textwidth]{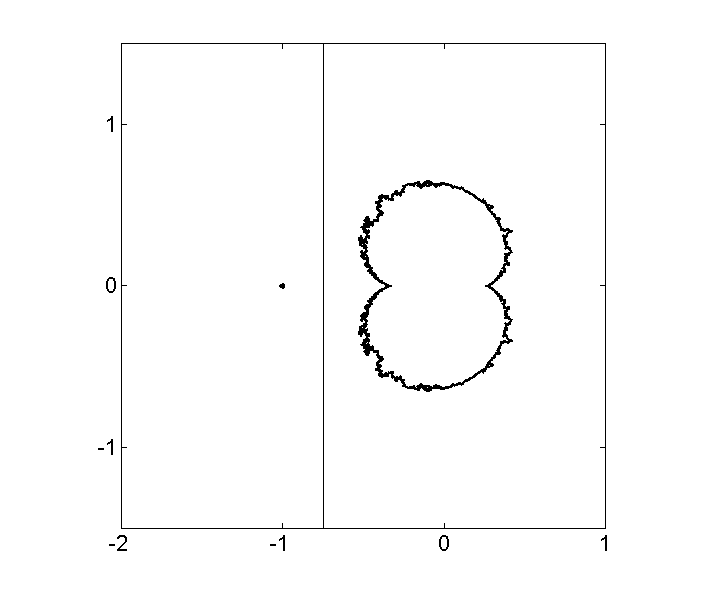}
\end{center}
\caption{\emph{\small {\bf Example of disconnected equi-M set for a network with $N=3$ nodes.} The curve traces the boundary of the equi-M set shown in Figure~\ref{connect}a, separated into two connected components by the line $Re(z) = -3/4$. The network is given by: $z_1 \to z_1^2+c$, $z_2 \to (az_1+z_2)^2+c$, $z_3 \to (z_1+z_2+bz_3)^2+c$, with connectivity weights $a = -1$, $b=-1$.}}
\label{connect2}
\end{figure}

\begin{prop}
The equi-M set for the network in the self-drive family above with connectivity weights $a = -1$, $b=-1$ is disconnected.
\label{example}
\end{prop}

\noindent {\bf Proof.} Notice first that, in general, all three node-wise projections of the critical orbits are real. We will show that the equi-M set described in the proposition has at least two connected components (the component of the origin and the component of $c=-1$), separated by the line $\text{Re}(z)=-3/4$. 

Indeed, the critical orbit is fixed for $c=0$, so that $c=0$ is trivially in the equi-M set of the network. Also, one can easily see that this particular self-drive network is postcritically finite when $c=-1$. Indeed, the first component of the critical orbit has in this case period two ($0 \to -1$); the second component has period four ($0 \to -1 \to -1 \to 0$) and the third component has period four ($0 \to -1 \to 0 \to 0$).

Finally, one can easily prove that no point on the line  $\text{Re}(c)=-3/4$ is in the equi-M set of this network. Indeed, notice that $c=-3/4$ is the only point in the traditional Mandelbrot set with $\text{Re}(c) =-3/4$ (it is the point joining the main cardioid and the period two bulb). Since the network M-set is a subset of the node-wise M-set for $z_1$ (which is the traditional Mandelbrot set), it also cannot contain any other points with $\text{Re}(c)=-3/4$. Furthermore, for our network, it can be shown that the third component $z_3$ of the critical orbit escapes when $c=-3/4$. Hence no point on the vertical line $c=-3/4$ is in the equi-M set of the network. It is interesting that the node that causes the pinch in the traditional M-set and renders the network M-set disconnected is in fact the ``output'' node, which receives control from both of the other two nodes; yet it is the orbit of $z_3$ that escapes, while the other two remain bounded when initiated at zero. Since the calculations are a little technical, we include them for completion in Appendix A.

\hfill $\Box$

\begin{figure}[h!]
\begin{center}
\includegraphics[width=0.5\textwidth]{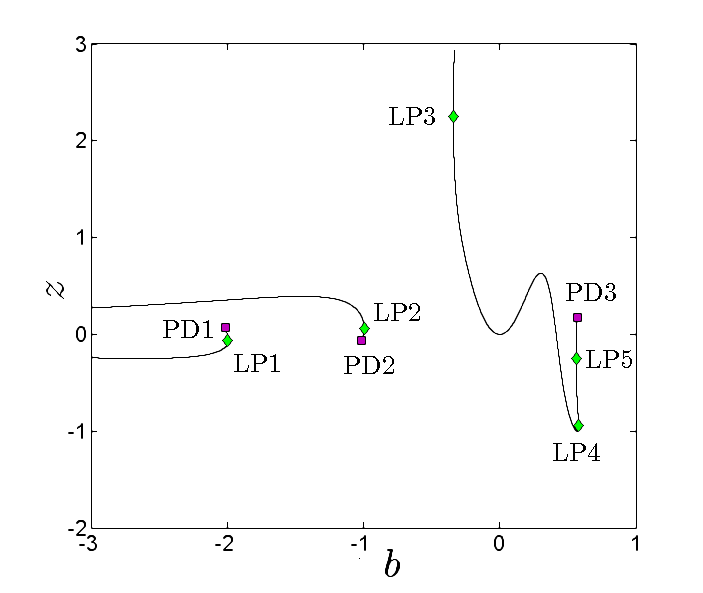}
\end{center}
\caption{\emph{\small {\bf Bifurcation diagram with respect to the coupling parameter $b$} for the function $f(\xi) = f_4 \circ f_3 \circ f_2 \circ f_1(\xi)$ which computes batches of four iterations of the node $z_3$, when $a=-1$, so that $f_1(\xi)=b^2\xi^2-1$, $f_2(\xi)=(b\xi-2)^2-1$, $f_3(\xi)=(b\xi-1)^2-1$ and $f_4(\xi)=f_3(\xi)$. The diagram shows three equilibrium curves, with a green diamond marking saddle node bifurcations (limit points/LP), and purple squares marking the first period doubling point of period doubling cascades to chaos. The intervals on which there is a stable equilibrium are: $b\sim-2.001$ (PD1) to $b\sim-1.995$ (LP1); $b\sim-1.01$ (PD2) to $b\sim-0.996$ (LP2); $b\sim-0.34$ (LP3) to $b\sim0.58$ (LP4), with a second stable fixed point between $b\sim0.56$ (LP5) and $b\sim0.57$ (PD3). The subsequent period doubling and chaotic windows are not shown, for clarity of the diagram (since two of these windows are extremely small), but the critical orbit remains bounded within this extended parameter range (as mentioned in the text).}}
\label{bifurcation_tree1}
\end{figure}

\vspace{5mm}
\noindent More generally, one can fix $c=-1$ and $a=-1$, keeping the critical orbits of the first two nodes periodic ($z_1$ performs a period two oscillation between $0 \to -1$, and $z_2$ has a period four oscillation $0 \to 0 \to -1 \to -1$), and study the effect of changing the self-drive parameter $b$ on the critical orbit of the node $z_3$. Recall that the trajectories of the critical orbits are real when $c$, $a$, $b$ are real. The bifurcation diagram in Figure~\ref{bifurcation_tree1} illustrates that in the parameter slice $a=-1$, there are at least three intervals for $b$ for which $c=-1$ is in the equi-M set. 

 These intervals include the windows corresponding to attracting fixed points, but also encompass period-doubling cascades and chaotic windows (not shown in Figure~\ref{bifurcation_tree1}). Given with three decimal approximation, these intervals for the parameter $b$ are: $[-2.016,-1.995]$, $[-1.028,-0.996]$ and $[-0.34,0.611]$.

\begin{figure}[h!]
\begin{center}
\includegraphics[width=0.7\textwidth]{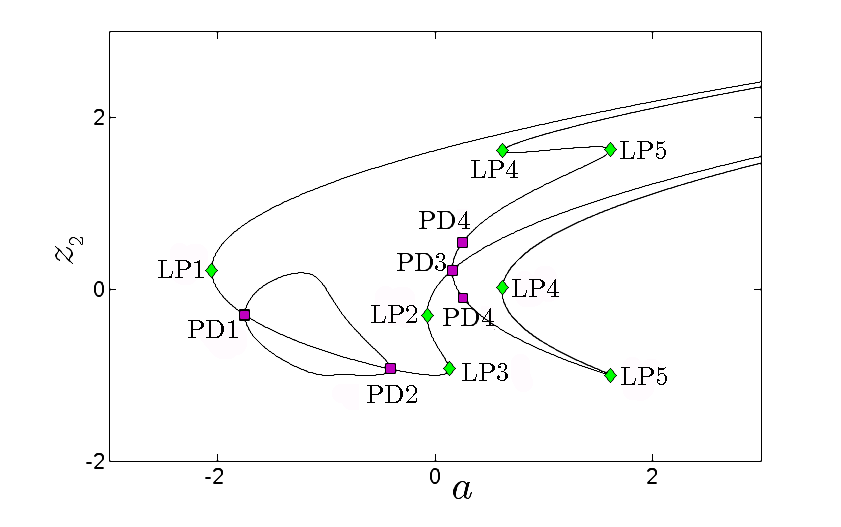}
\end{center}
\caption{\emph{\small {\bf Bifurcation diagram with respect to the coupling parameter $a$} for the function $f(\xi) = (\xi^2-1+a)^2-1$ (representing the even iterations of the $z_2$ component of the critical orbit in the feed-forward family when $c=-1$.}}
\label{bifurcation_tree2}
\end{figure}

\vspace{5mm}
\noindent Even more generally, one can compute the range for the coupling parameter $a$ which guarantees that $c=-1$ remains in the node-wise Mandelbrot set for the second node $z_2$. Since the critical orbit is real for real values of $c$ and $a$, we can study this by tracking the bifurcations of the function $f(\xi) = (\xi^2-1+a)^2-1$ (which represents the transition between even iterations, as calculated above) with respect to $a$. We show the bifurcation diagram schematically in Figure~\ref{bifurcation_tree2}. The initial condition $z_2(0)=0$ escapes for $a<-2$; it converges to a stable fixed point, and then to a stable periodic two orbit (after the period doubling at $a=-5/4$). The attracting period two orbit survives (and attracts the origin) until $a=-0.4$ (with a superattracting stage at $a=-1$), and then collapses back into a stable fixed point. As $a$ is increased, the system undergoes a cascade of period doubling bifurcations, starting with the first one at $a~\sim 0.15$, birthing periodic cycles which continue to attract the critical orbit; it continues along on root to chaos, maintaining $z_2$ bounded within $[-1,2]$. Eventually, the origin escapes the trapping interval when the parameter crosses the value $a~\sim 0.7$. 

\begin{figure}[h!]
\begin{center}
\includegraphics[width=\textwidth]{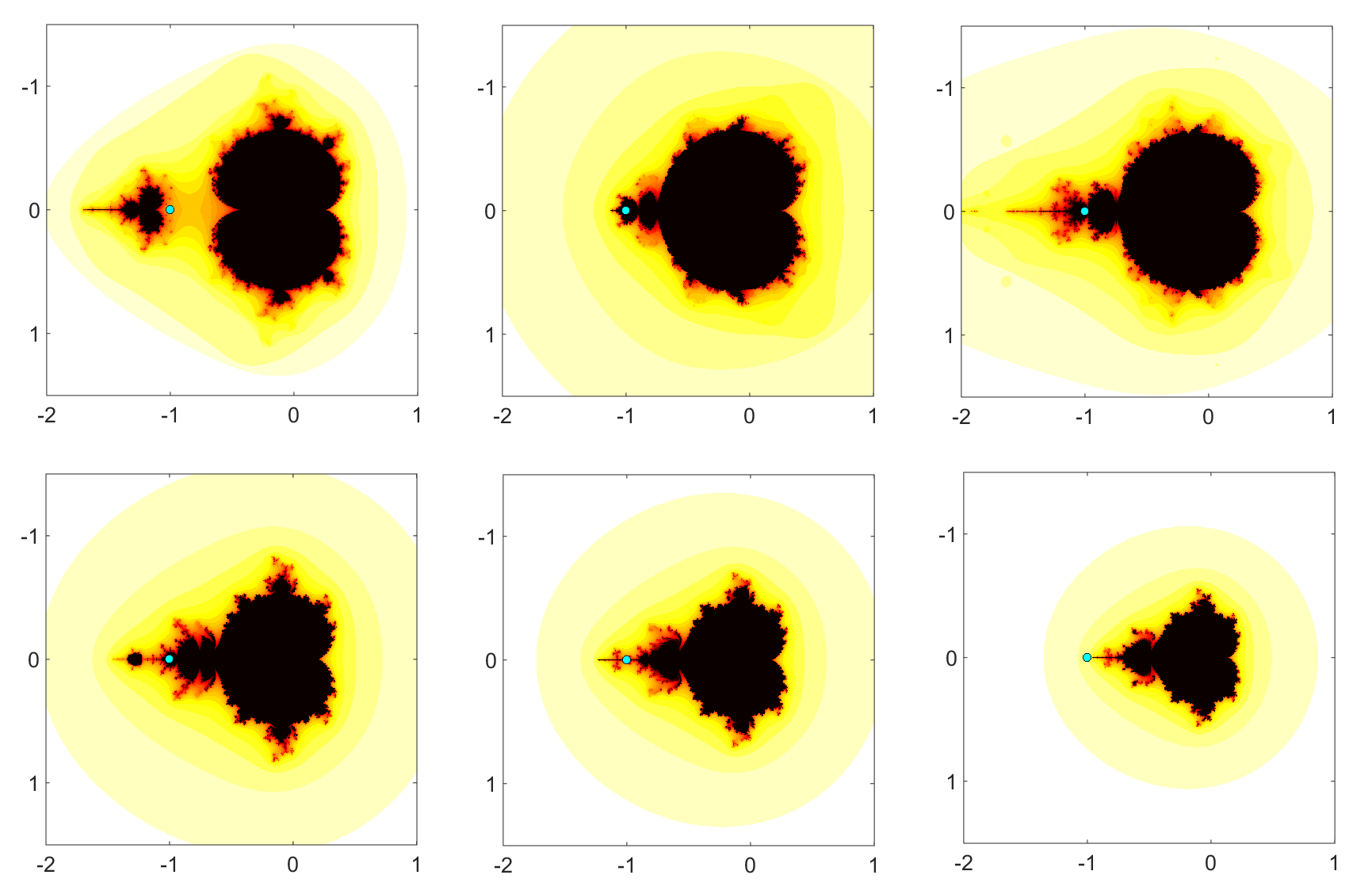}
\end{center}
\caption{\emph{\small {\bf Node-wise Mandelbrot sets for $z_2$}, illustrated for different values of $a$ along the bifurcation diagram in Figure~\ref{bifurcation_tree1}. In each equi-M set, the cyan dot represents the point $c=-1$. {\bf Top.} From left to right: $a=-2.1$ (M-set is pinched at $c=-1$ and $z_2$ component of the critical point escapes); $a=-1$ (super-attracting orbit of period two at $c=-1$) ; $a=-0.5$ ($z_2$ component of critical point is pre-periodic at $c=-1$). {\bf Bottom.} $a \sim 0.2$ (super-attractive orbit of period two at $c=-1$); $a~\sim 0.4$ (super-attractive orbit of period five at $c=-1$);  $a=0.75$ (M-set falls short of $c=-1$ and $z_2$ component of the critical point escapes).}}
\label{bif_values}
\end{figure}

Hence $c=-1$ is in the $z_2$ Mandelbrot set for the relatively large interval $[-2,0.7]$ for $a$, and is not in the node-wise Mandelbrot set outside of this parameter range. The two endpoints of this interval have different significance and mechanisms. On one hand, when lowering $a$ past the low critical state $a=-2$, the point $c=-1$ pinches out and separates the $z_2$ Mandelbrot set into two connected components (to the right and to the left of the line Re$(z)=-1$, see Figure~\ref{bif_values}a). On the other hand, when raising $a$ in the positive range, the tail of the $z_2$ Mandelbrot set shortens, so that past the high critical state $a=0.7$, the point $c=-1$ is left out, and the whole set is to the right of the line Re$(z)=-1$ (see Figure~\ref{bif_values}b). Along this interval, there are values of $a$ for which the node $z_2$ has a super-attracting orbit at $c=-1$ (for example, the critical orbit is periodic at $a=-2$, $a=-1.6$, $a=-1$, $a=0$, $a=1/5$, $a=0.3$). There are also values of $a$ for which the $z_2$ component of the critical orbit is pre-periodic at $c=-1$ ($a=-0.5$).

One can then further look at the third component of the critical orbit corresponding to the node $z_3$. For example, in the case when the critical orbit of $z_2$ stabilizes asymptotically to an attracting period two oscillation, this oscillation is represented by a fixed point $\xi_0$ for the function $f(\xi)=(\xi^2-1+a)^2-1$ above. Hence four batch iterations of $z_3$ converge to the function $g(\xi) = (\xi_0+b\xi)^2-1$. For $\xi_0$ in the intervals found above, one can study asymptotic dependence on $b$, by constructing bifurcation diagrams similarly to that in Figures~\ref{bif_values}.

\begin{figure}[h!]
\begin{center}
\includegraphics[width=\textwidth]{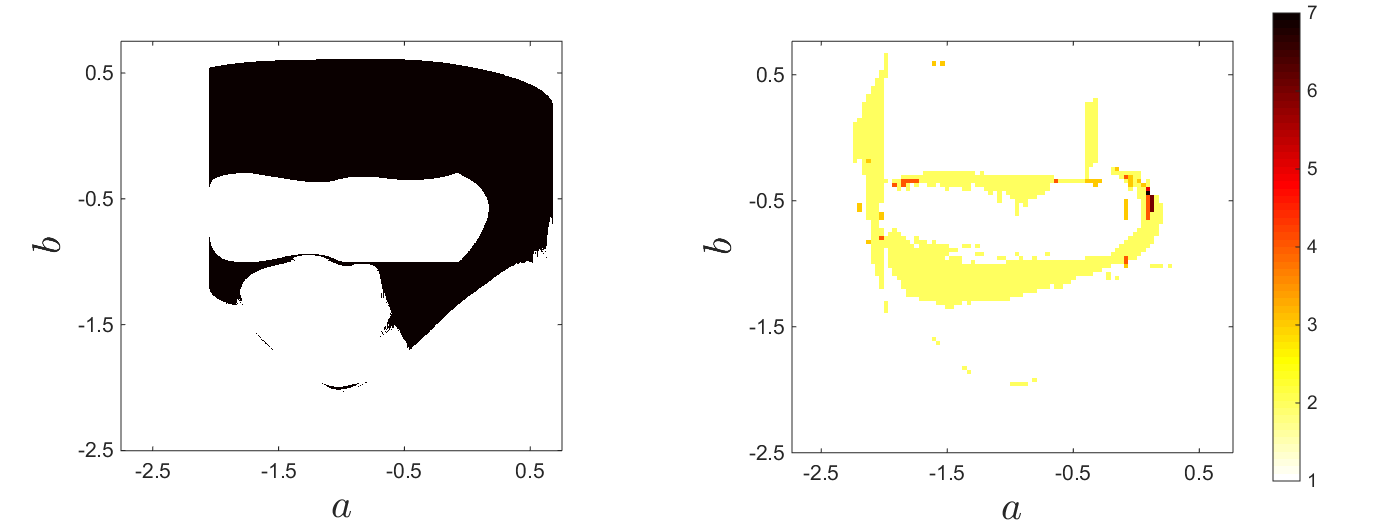}
\end{center}
\caption{\emph{\small {Parameter loci in the $(a,b)$ plane,} for the network family given by: $z_1 \to z_1^2+c$, $z_2 \to (az_1+z_2)^2+c$, $z_3 \to (z_1+z_2+bz_3)^2+c$. {\bf Left.} Locus (computed within the rectangle $[-2.75,0.75] \times [-2.75,0.75]$, shown in black)  of pairs $(a,b)$ for which $c=-1$ is in the equi M-set. {\bf Right.} Connectivity locus of the equi-M set, within the same rectangle $[-2.75,0.75] \times [-2.5,0.75]$, computed using the blowup algorithm before assessing connectivity of the sets. The color corresponding to each $c$ represents the estimated number of connected components of the equi-M set (as shown in the color bar).}}
\label{locus}
\end{figure}

\vspace{5mm}
\noindent Based on these observations, one may investigate if there is a relationship between how the centers of the former hyperbolic components  of the Mandelbrot set are being perturbed by the network structure (and whether they still belong to the equi-M set) and the connectedness of the equi-M set as a whole.  Below, we try to understand this comparison, using a more comprehensive illustration of asymptotic behavior within the particular three-dimensional family of feed-forward networks considered above. In Figure~\ref{locus}a, we show the connectivity parameter locus $(a,b)$ (represented along the horizontal and respectively vertical coordinate axes) for which the complex parameter $c=-1$ is in the equi-M set. In Figure~\ref{locus}b, we show the result of a rough computation of the number of connected components in the equi-M set, for each parameter pair $(a,b)$.  Due to difficulty in the reduced resolution (that was necessary for insuring feasible computation time) we used a ``blow-up'' algorithm that expanded each equi-M set by small margin before assessing its connectedness. While this may be introducing some negative error in detecting distinct connected components, we found that it substantially reduces positive detection error (due to the inability of the numerical code to identify filaments in the original equi-M sets represented in reduced resolution, as further explained in Appendix B).

It is interesting to reinterpret the bifurcation diagram in Figure~\ref{bifurcation_tree1} in the broader context of Figure~\ref{locus}a. The former represents the slice $a=-1$ of the latter, so that one can observe the three black intervals for $b$ along the vertical line $a=-1$, representing the three windows in the bifurcation diagram where the critical orbit is bounded for $c=-1$. It is also interesting, although less trivial, to compare the left and right panels of Figure~\ref{locus}. Although the presence of $c=-1$ in the equi-M set does not imply connectedness, it is clearly related to the connectedness locus, with a break in connectedness (around the yellow region representing the boundary between one and two connected components in Figure~\ref{locus}a) seemingly related to the boundary of the inner white region in Figure~\ref{locus}a, where $c=-1$ is pinched out of the M-set).

\subsection{Fatou-Julia Theorem extended}

\begin{figure}[h!]
\begin{center}
\includegraphics[width=\textwidth]{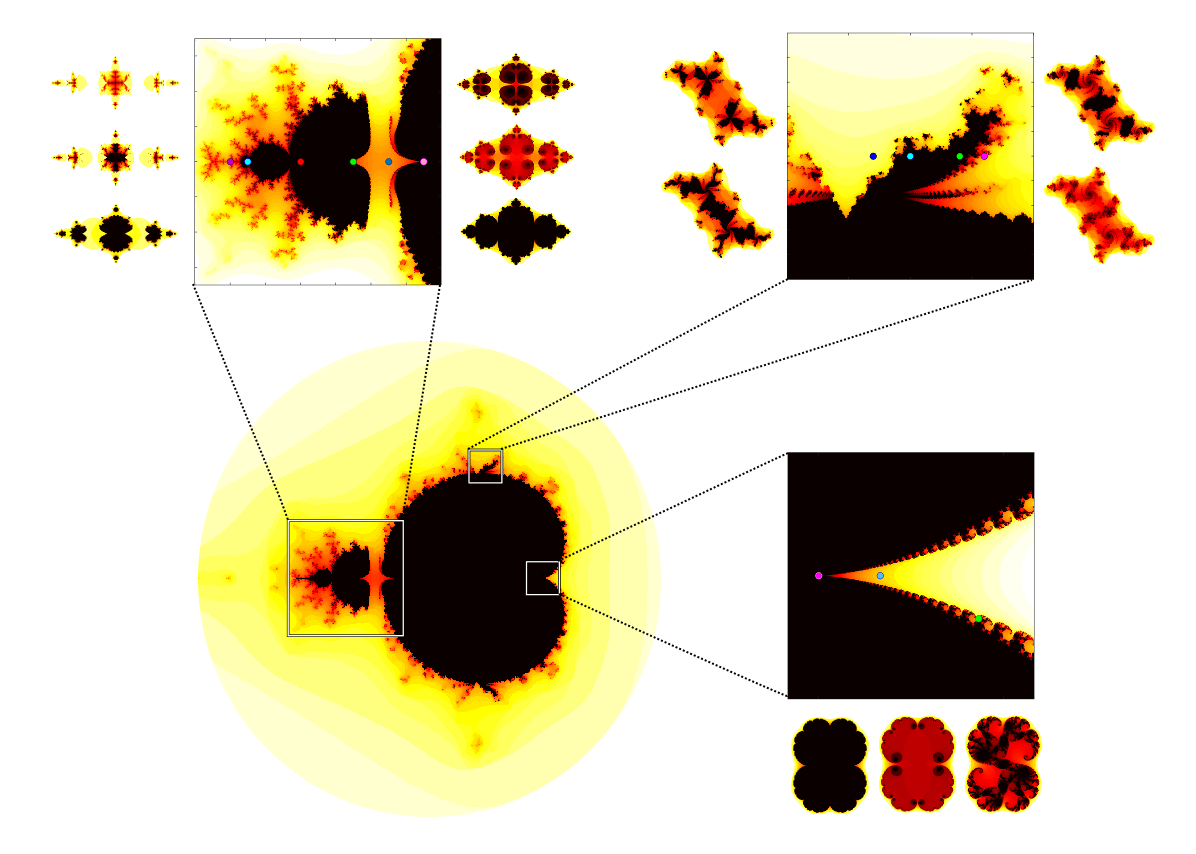}
\end{center}
\caption{\small \emph{{\bf Uni-Julia sets} for a self-drive network with $a=-1/3$ and $b=-1/3$, for different values of the equi-parameter $c$. We magnified three rectangular windows around the boundary of the network's equi-M set: $[-0.2,0.2] \times [0.245,0.285]$ (around the cusp), $[-0.2,0] \times [0.58,0.78]$ (top) and $[-1.3,-0.6] \times [-0.3,0.3]$ (around the tail). For each window, we show several uni-Julia sets corresponding to the $c$ values marked in colors. For each magnification window, as the dots are listed from left to right, the corresponding uni-Julia sets are represented from left to right and then top to bottom.}}
\label{connect_illustr}
\end{figure}

In our previous work, we noticed that existence of uni-Julia sets with finitely many connected components breaks, in the case of networks, the connected/dust duality on which the Fatou-Julia theorem is based in the traditional case of single iterated maps. In the same reference, we relied on a few numerical  illustrations of uni-Julia sets for a variety of parameters $c$ (see for example Figure~\ref{connect_illustr}), chosen both inside and outside of the equi-M set for their respective network, to conjecture that the uni-J set is connected only if $c$ is in the equi-M set of the network, and it is totally disconnected only if $c$ is not in the equi-M set of the network. 


Here we illustrate this relationship in greater detail, while still using numerical approaches. In Figure~\ref{connect_illustr}, we show the equi-M set for one of our self-drive example networks, together with a the uni-J sets corresponding to a collection of points $c$ chosen close to the boundary of the equi-M set (so that some of them are inside the equi-M set, and some are outside). The illustration supports the idea that, although the connectivity of the uni-J sets (from one, to finitely many, to infinitely many components) degrades in the proximity of the boundary of the M-set, there is no sudden break that happens precisely on the boundary, like in the case of single map iterations.

For a more systematic view, we computed and illustrated together, for a few example networks, the boundary of the equi-M set and the connectedness locus for the uni-J set. While the former was relatively easy to compute as the critically bounded locus for the network, the latter presented some difficulties in reconciling computational efficiency with obtaining uni-J sets in sufficiently good resolution to allow us to estimate their number of connected components. This was problematic in particular for the situations where the Julia set had short, thin filaments, likely to escape detection in low resolution, in which case we suspected the code to report ``fake'' connected components, and thus over-count the number of components. To eliminate this positive error without increasing the $z$-plane resolution (which impacts computational time quadratically), we used a common ``blow-up'' technique, adding a small border to each uni-J set to account for the possible connections due to filaments. This, of course, may introduce the opposite type of error (that of under-counting components). However, the two methods produced unexpectedly similar results qualitatively, in the sense of identifying the same loci of connectedness and total disconnectedness. In the transitional region, the connected component counts were higher with the first algorithm versus the second, as one would have expected (see Appendix B). 

\begin{figure}[h!]
\begin{center}
\includegraphics[width=\textwidth]{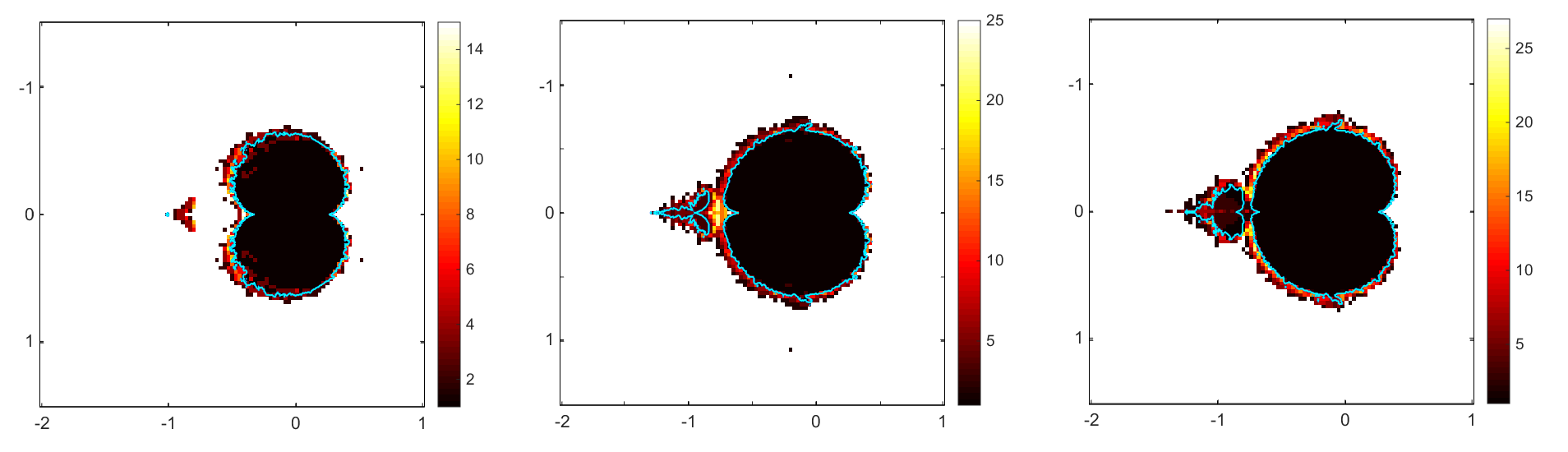}
\end{center}
\caption{\small \emph{{\bf Comparison between the equi-M set and the uni-J connectedness locus} for the self-drive networks illustrated in Figure~\ref{connect}: {\bf A.} $a=-1$, $b=-1$; {\bf B.} $a=-1/3$, $b=-1/3$; {\bf C.} $a=-2/3$, $b=-1/3$. The panels represent the square $[-2,1] \times [-1.5,1.5]$ in the equi-parameter plane. The cyan curve represents the boundary of the equi-M set, computed with 50 iterations. The colors correspond to the number of connected components for the respective uni-J set (computed approximately using the numerical algorithm discussed in Appendix B), with the color scheme going from white (inside region, one connected component) through tones of red and yellow, as the number of finitely many connected components increases to 2, 3, etc (see color bar). White corresponds to the locus where the uni-J set was found to be dust (the numerical could not capture the totally disconnected points, so it returned the answer as ``zero'' components, which we then scaled by hand to appear as white background).}}
\label{connect_compare}
\end{figure}

In Figure~\ref{connect_compare}, we illustrate the implementation of the blow-up algorithm on the three self-drive networks shows in Figure~\ref{connect}. We computed both the equi-M set in $\mathbb{C}$ (in the sense defined in Section~\ref{prior}), as well as the connectedness locus (also in $\mathbb{C}$) for the uni-J set of the network. While it does not come as a surprise that the two are no longer identical, we found that they are clearly related. Future work will focus on obtaining an analytic understanding of this relationship.


\section{Network methods}
\label{networks}
\subsection{Spectral versus dynamic classes}
\label{sets}

For small or very simple networks, one can try to identify specifically the effect of different graph architectural properties onto the ensemble asymptotic dynamics. As we have done in previous work for continuous time systems, we first investigate possible relationships between the network adjacency spectrum and the class of ensemble dynamics. For a network with discrete quadratic nodes, it seems natural to characterize the network by the properties of its asymptotic sets: the equi-M set, and the uni-Julia set, for a fixed equi-parameter $c$.
 
\begin{defn}
We say that two networks ${\cal N}_1$ and ${\cal N}_2$ are in the same asymptotic class if, for any initial condition $(z_0^1,z_0^2, \hdots, z_0^N) \in \mathbb{C}^N$, its multi-orbit under ${\cal N}_1$ iterates out of the escape disc at the same rate as when iterated under ${\cal N}_2$. We say that they are in the same uni-asymptotic class if the same applies for the multi-orbits of all uni-initial conditions $z_0 \in \mathbb{C}$.
\end{defn}

\noindent {\bf Remark.}  Visually, this means that the corresponding prisoner sets (or uni-prisoner sets, respectively) are identical between two networks in the same uni-asymptotic class, and so are the escape sets, with identical ``escape colors'' assigned to corresponding points.

\begin{conj}
Network uni-asymptotic classes are invariant under changes of the equi-parameter.
\end{conj}

\noindent The conjecture states a potentially very useful result: that two distinct configurations which produce identical uni-asymptotic dynamics for one value of the parameter $c$, will also do so for all any other value of $c$, and two configurations which produce different asymptotic structure under one value of $c$ will still do so under any other value of $c$.

To investigate this hypothesis numerically, we focused on replicating the result in networks with two types of general restrictions: (1) networks with a fixed number of nodes $N$ and a fixed number of edges $j$, with no additional conditions on the configuration; (2) bipartite networks with $N$ nodes in each of the two interconnected cliques (previously used to represent interacting neural populations in our modeling work), and with specified number of edges $i$ and $j$ between the two cliques, respectively. In Appendix C, we illustrate one example from each category. 

\begin{figure}[h!]
\label{table_1_3}
\begin{footnotesize}
\begin{tabular}{|c|c|c|c|}
\hline
& & &\\
$\;\left[ \begin{array}{cc|cc} &  & 1 & 0\\  &  & 0 & 0 \\ \cline{1-4} 1 & 1 &  & \\  1 & 0 &  & \end{array} \right] \; ({\cal A}_i)$ &
$\;\left[ \begin{array}{cc|cc} &  & 0 & 1\\  &  & 0 & 0 \\ \cline{1-4} 1 & 1 &  & \\  1 & 0 &  & \end{array} \right] \; ({\cal B}_{ii})$ &
$\;\left[ \begin{array}{cc|cc} &  & 0 & 0\\  &  & 1 & 0 \\ \cline{1-4} 1 & 1 &  & \\  1 & 0 &  & \end{array} \right] \; ({\cal B}_{iii})$ &
$\;\left[ \begin{array}{cc|cc} &  & 0 & 0\\  &  & 0 & 1 \\ \cline{1-4} 1 & 1 &  & \\  1 & 0 &  & \end{array} \right] \; ({\cal C}_{iv})$\\
& & &\\
\hline
& & &\\
$\;\left[ \begin{array}{cc|cc} & & 1 & 0\\  &  & 0 & 0 \\ \cline{1-4} 1 & 1 &  & \\  0 & 1 &  & \end{array} \right] \; ({\cal B}_{iii})$ &
$\;\left[ \begin{array}{cc|cc} & & 0 & 1\\  &  & 0 & 0 \\ \cline{1-4} 1 & 1 &  & \\  0 & 1 &  & \end{array} \right] \; ({\cal C}_{iv})$ &
$\;\left[ \begin{array}{cc|cc} & & 0 & 0\\  &  & 1 & 0 \\ \cline{1-4} 1 & 1 &  & \\  0 & 1 &  & \end{array} \right] \; ({\cal A}_i)$ &
$\;\left[ \begin{array}{cc|cc} & & 0 & 0\\  &  & 0 & 1 \\ \cline{1-4} 1 & 1 &  & \\  0 & 1 &  & \end{array} \right] \; ({\cal B}_{ii})$ \\
& & &\\
\hline
& & &\\
$\;\left[ \begin{array}{cc|cc} & & 1 & 0\\  &  & 0 & 0 \\ \cline{1-4} 1 & 0 &  & \\  1 & 1 &  & \end{array} \right] \; ({\cal B}_{ii})$ &
$\;\left[ \begin{array}{cc|cc} & & 0 & 1\\  &  & 0 & 0 \\ \cline{1-4} 1 & 0 &  & \\  1 & 1 &  & \end{array} \right] \; ({\cal A}_i)$ &
$\;\left[ \begin{array}{cc|cc} & & 0 & 0\\  &  & 1 & 0 \\ \cline{1-4} 1 & 0 &  & \\  1 & 1 &  & \end{array} \right] \; ({\cal C}_{iv})$ &
$\;\left[ \begin{array}{cc|cc} & & 0 & 0\\  &  & 0 & 1 \\ \cline{1-4} 1 & 0 &  & \\  1 & 1 &  & \end{array} \right] \; ({\cal B}_{iii})$ \\
& & &\\
\hline
& & &\\
$\;\left[ \begin{array}{cc|cc} & & 1 & 0\\  &  & 0 & 0 \\ \cline{1-4} 0 & 1 &  & \\  1 & 1 &  & \end{array} \right] \; ({\cal C}_{iv})$ &
$\;\left[ \begin{array}{cc|cc} & & 0 & 1\\  &  & 0 & 0 \\ \cline{1-4} 0 & 1 &  & \\  1 & 1 &  & \end{array} \right] \; ({\cal B}_{iii})$ &
$\;\left[ \begin{array}{cc|cc} & & 0 & 0\\  &  & 1 & 0 \\ \cline{1-4} 0 & 1 &  & \\  1 & 1 &  & \end{array} \right] \; ({\cal B}_{ii})$ &
$\;\left[ \begin{array}{cc|cc} & & 0 & 0\\  &  & 0 & 1 \\ \cline{1-4} 0 & 1 &  & \\  1 & 1 &  & \end{array} \right] \; ({\cal A}_i)$ \\
& & &\\
\hline
\end{tabular}
\end{footnotesize}

\vspace{3mm}
\hspace{-2mm}
\includegraphics[width=0.96\textwidth]{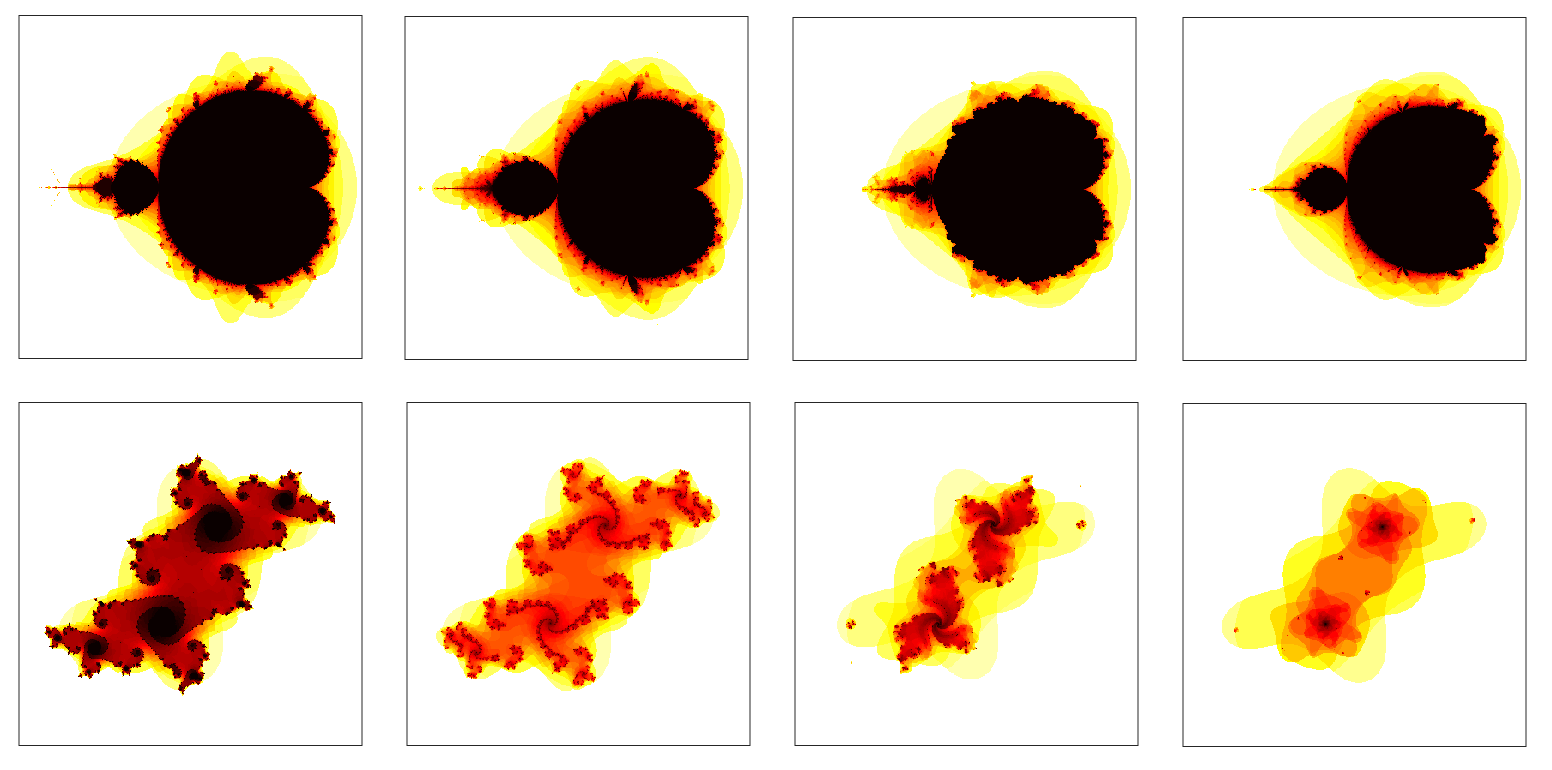}

\caption{\label{table_3_3} \emph{\small{\bf Adjacency and dynamics classes for N=2, density type ($M_{xy}$,$M_{yx}$)=(1,3) and $g_{xx}=g_{yy}=0.5$, $g_{xy}=g_{yx}=-0.5$.}  Adjacency classes are designated by letters (${\cal A} - {\cal C}$) and asymptotic classes denoted by the subscript ($i - iv$). The top figure panels represent the equi-M sets for all asymptotic classes $i$-$iv$. The bottom figure panels show the $i$-$iv$ uni-planes for the equi-parameter ($c=-0.117-0.856i$), with prisoners plotted in black and escapees plotted in colors according to the escape rate. Notice that in this case one can achieve all dynamics classes by changing either one of the diagonal block matrices, while keeping the other fixed.}}
\label{classes_bipartite}
\end{figure}

In Figure~\ref{classes_simple}, we considered all networks with $N=3$ nodes and $j=7$ edges, with all edge weights set as $g=1/N$. The panels illustrate the uni-J sets for the equi-parameter values $c=-1.15+0.26i$ and $c=-0.13+i$. In Figure~\ref{classes_bipartite}, we considered all bipartite networks with $N=2$ nodes per clique, $i=1$ and $j=3$, with positive weights $g=1/2$ for the edges connecting nodes within the cliques, and negative weights $g=-1/2$ for the edges between the cliques. The panels illustrate the uni-J and equi-M sets for the equi-parameter value $c=-0.117-0.856i$.

Spectral and asymptotic classes are not in a one-to-one correspondence, either way. Notice, in both tables, that two distinct matrices from the same spectral class may produce in some cases identical, in other cases different uni-asymptotic dynamics. Conversely, two matrices in different spectral classes may produce the same uni-asymptotic dynamics. However, even though not determined by the adjacency spectrum, uni-asymptotic classes remain consistent for all value of $c$.

\subsection{Core asymptotic sets}
\label{average}

\noindent In previous work, we have explored a statistical approach to relating graph structure to asymptotic dynamics in networks~\cite{radulescu2015nonlinear}. When interested in all network configurations with a specific property ${\cal P}$ (e.g., density of oriented edges), one may consider, for each initial point (or alternately for each point in parameter space) the fraction of all configurations which produce a specific asymptotic behavior. Then, a ``probabilistic'' bifurcation can be defined in terms of the likelihood of a system to transition  between two different behaviors when the edge configuration is slightly perturbed, when the only knowledge we have on the network configuration is property ${\cal P}$.    

For example, fix an equi-parameter $c$, and, for simplicity, set all edge weights in the network equl to $1/N$, where $N$ is the size of the network. Consider the property ${\cal P}$ to be fixing the number of edges to $k$, with $0 \leq k \leq N^2$. For each $z_0 \in \mathbb{C}$, we count the fraction of configurations with ${\cal P}$ for which the uni-orbit of $z_0$ is bounded. 

\begin{defn} We call the \textbf{core uni-prisoner set} the set of all points $z_0 \in \mathbb{C}$, for which the initial condition $(z_0,...z_0) \in \mathbb{C}^N$ produces a bounded multi-orbit when iterated under all network configurations with property ${\cal P}$. We call the \textbf{core uni-J set}\footnote{This term was chosen in order to emphasize the analogy with a similar concept defined by Sumi in the case of random iterations of postcritically bounded polynomials~\cite{sumi2006semi,sumi2015random}} the boundary of this set in $\mathbb{C}$.
\end{defn}

\noindent Instead of inspecting connectivity of each configuration-specific uni-J set at a time, one can instead study topological properties of the level sets of ${\cal P}$ in the complex $z$-plane, in particular connectivity of the core uni-J set (which is the boundary of the $1$-level set). One can track how the core uni-J set is affected when changing the edge weights $g$, the equi-parameter $c$, the network size $N$ or, finally, even the network fixed property ${\cal P}$. Furthermore, one can distinguish between the parameter values for which the core uni-Julia set remains connected for all edge configurations with property ${\cal P}$, versus parameter values for which changes in edge density alter connectivity of the core uni-J set.\\

\begin{figure}[h!]
\begin{center}
\includegraphics[width=0.7\textwidth]{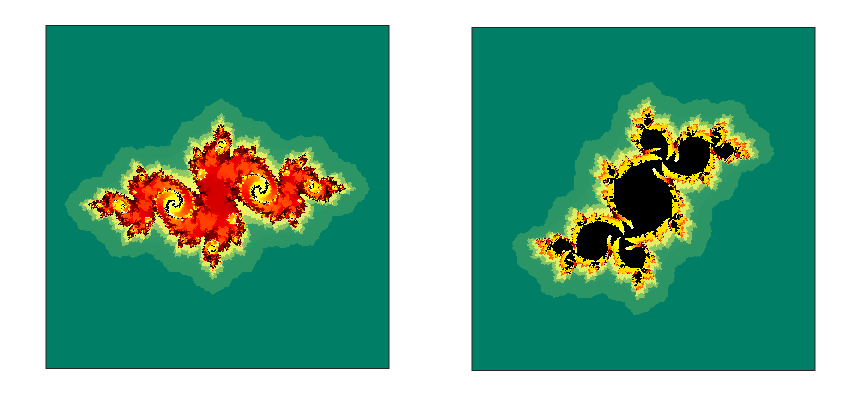}
\end{center}
\caption{\emph{\small {\bf Core uni-J sets over all network configurations with $N=3$ nodes, edge density $\delta=7/9$, for fixed edge weights $g=1/3$}, and fixed equi-parameter $c$. {\bf A.} $c=-1.15+0.26i$; {\bf B.} $c=-0.13+i$. All panels were computed for $50$ iterations, with spacial resolution $200 \times 200$, and escape radius $R_e=20$.}}
\label{prob_Julia_N3_j_7}
\end{figure}

\begin{figure}[h!]
\begin{center}
\includegraphics[width=0.7\textwidth]{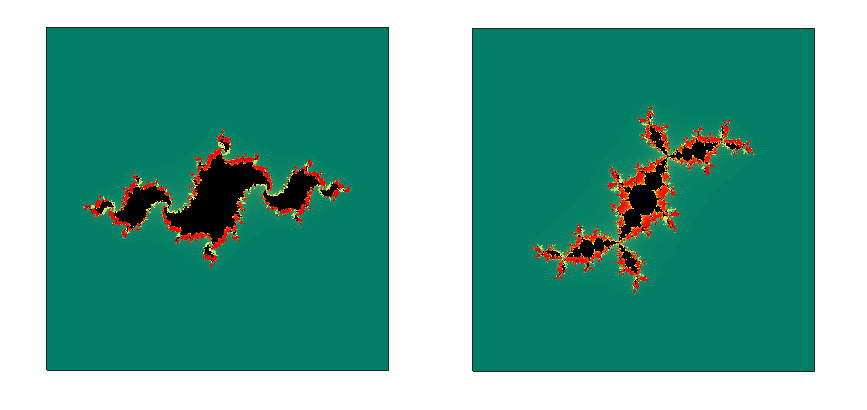}
\end{center}
\caption{\emph{\small {\bf Core uni-J sets over all network configurations with $N=3$ nodes, edge density $\delta=8/9$, for fixed edge weights $g=1/3$}, and fixed equi-parameter $c$. {\bf A.} $c=-1.15+0.26i$; {\bf B.} $c=-0.13+i$. All panels were computed for $50$ iterations, with spacial resolution $200 \times 200$, and escape radius $R_e=20$.}}
\label{prob_Julia_N3_j_8}
\end{figure}

\noindent To fix out ideas, we discuss the concept of core uni-J set in the case of property ${\cal P}$ being ``fixed edge density $\delta = k/N^2$.'' Figures~\ref{prob_Julia_N3_j_7} and~\ref{prob_Julia_N3_j_8}  illustrate core uni-J sets in networks of size $N=3$ and uniform edge weights $g=1/3$, for different equi-parameters $c$,  and different edge densities $\delta$. The color associated to each point $z_0 \in \mathbb{C}$ represents the likelihood (over all network configurations) for the initial condition $(z_0,z_0,z_0)$ to remain bounded under iterations of a network with node-wise dynamic specified by $c$ and edge density specified by $\delta$. In particular, the black central region represents the core uni-prisoner set. For example, Figure~\ref{prob_Julia_N3_j_7}a and b show the core Julia sets corresponding to the two classes of asymptotic dynamics described respectively in the left and right columns of Figure~\ref{classes_simple} in Appendix C.

Intuitively speaking, as one would expect, the network dynamics becomes generally more rigid for higher edge densities $\delta$, and more fluid for lower densities, since more edges are expected to increase communication and ``synchronization'' between nodes. This effect is clearly captured in the comparison between the lower density $\delta=7/9$ panels in Figure~\ref{prob_Julia_N3_j_7}, and the corresponding panels for the next higher density $\delta=8/9$ in Figure~\ref{prob_Julia_N3_j_7}. The level curves appear closer together in the higher density sets, so that small perturbations in the initial condition can more dramatically change the likelihood of a multi-orbit to escape. At a finer level, however, one can clearly notice that the effect of increasing the edge density $\delta$ on the core uni-J set varies with the network. For example, depending on the equi-parameter $c$, the core uni-prisoner set may gain in area and connectedness with increasing edge density (as seen in left panels of the two figures), or may shrink (as in the right panels). One can define and investigate the same concept similarly in equi-M sets:

\begin{defn} We call the \textbf{core equi-M set} the set of all points $c \in \mathbb{C}$ for which the critical multi-orbit is bounded in $\mathbb{C}^n$, when computed for all network configurations with property ${\cal P}$.
\end{defn}

For small networks, the equi-M set is highly sensitive to small changes in the network architecture, as one can see for example the Appendix C illustrations. By simply adding, deleting or moving one single edge, one can transition between asymptotic classes, thus altering substantially the geometry and properties of the equi-M set, and of the uni-J sets for all values of the parameter $c$. One is interested to ask the same type of questions in the context of higher-dimensional networks. Do small perturbations in the architecture affect the asymptotic behavior to a similar extent, or do the rest of the edges stabilize the network? Does the presence of  this``vulnerability'' depend on global properties such as overall edge density, or on local information, such as on the place where the addition/removal happened? These are important theoretical questions which relate to counterparts in modeling and the  life sciences.

\begin{figure}[h!]
\begin{center}
\includegraphics[width=0.8\textwidth]{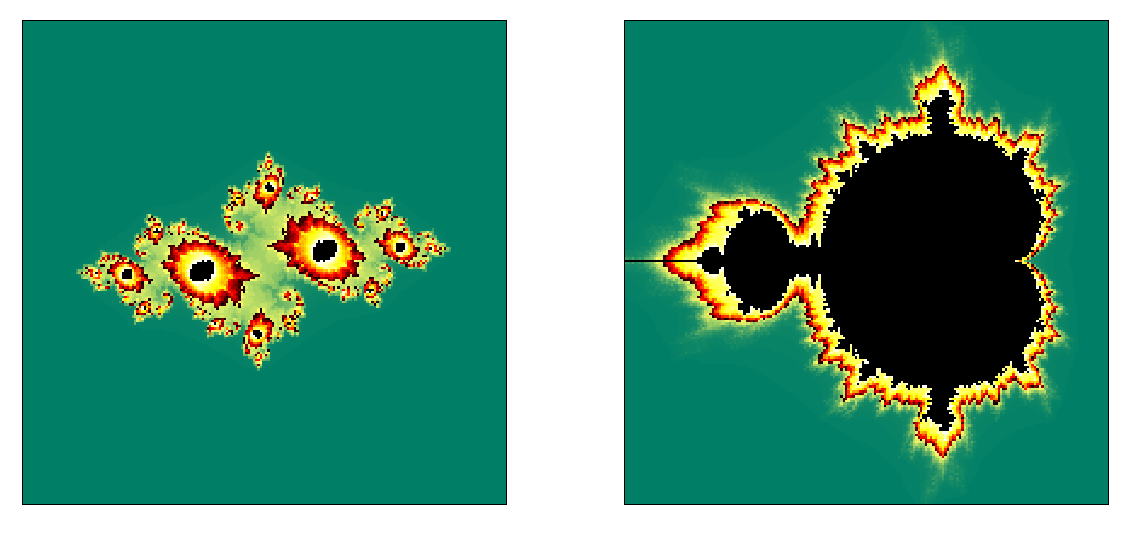}
\end{center}
\caption{\emph{\small {\bf Core sets for network configurations with $N=10$ nodes.} {\bf Left.} Core uni-J set for fixed equi-parameter $c=-1.15+0.26i$, edge density $\delta=80/100$, and fixed edge weights $g=1/N$. {\bf  Right.} Core equi-M set for edge density $\delta=60/100$ and edge weights $g=1/N$. All panels were computed for $k=50$ iterations, with spacial resolution $200 \times 200$, and escape radius $R_e=20$.}}
\label{high_dim}
\end{figure}

In Figure~\ref{high_dim}, we show a core uni-J set and the core equi-M set for the collection of all networks of $N=10$ nodes, with ${\cal P}$ being that they common edge density  $\delta= 80/100$. The total number of configurations with property ${\cal P}$ is extremely large (for the network size $N=10$, which is still relatively small, one obtains $\ds {100 \choose 80}$, which is of the order $10^{20}$). Even considering the equivalence classes of asymptotic dynamics (assuming we have identified them and their size), averaging over all possibilities is extremely challenging computationally. In our previous work, we have shown that sample-based means are quite accurate, even for very small samples. In Figure~\ref{high_dim}, we used samples of size ${\cal S}=20$ configurations out of the total of approximately $5 \times 10^{20}$ to illustrate our core sets.

These types of illustrations offer concomitant (while sample-based) stochastic information on the asymptotic dynamics within a large collection of networks. They could be important in that they may help detect asymptotic properties which are robust to changes in architecture, and distinguish them from those which are sensitive to change. For example, one of the features which we had previously noticed consistently is the persistence of the cups structure in all equi-M sets (with small variations in its position, depending on network architecture and node-wise dynamics). Figure~\ref{high_dim}b confirms this observation, showing minimal variability in the cusp area compared to regions of high sensitivity (such as the tail area, where even a small change in $c$ may lead from certainty (black) to very small likelihood (yellow and red)  of a bounded critical orbit.


\section{Discussion}
\label{discussion}

\subsection{Specific comments}

In this paper, we reformulated some well-known questions from single map quadratic dynamics in the context of iterations of  ensemble quadratic maps, coupled up in a network, according to an underlying adjacency graph structure. We investigated whether single map results regarding orbit convergence, escape radius and the topological structure of asymptotic sets change when studying a small network of $n$ quadratic complex nodes. We focused in particular on one-dimensional complex slices of these sets in $\mathbb{C}^n$, which we call uni-J sets and equi-M sets. 

We found that, while some of the structure of the traditional Mandelbrot set is conserved in Mandelbrot slices (such as fractality on the boundary, or the cusp at its rightmost point along the real axis) the equi-M sets no longer exhibit the hyperbolic bulb structure, and are no longer necessarily connected. In fact, depending on the architecture of the network and the strength of the connections between nodes, the original centers of the hyperbolic components may no longer be within the equi-M set altogether.

Similarly, the connection between the Julia and Mandelbrot set is a lot more complicated in networks of nodes. We have investigated a variant of the equivalence between the Mandelbrot set and the connectedness locus of the Julia set, as originally stated for single maps. Since relating the network Julia and Mandelbrot sets as loci in $\mathbb{C}^n$ seemed rather difficult, we started by comparing the structure of the equi-M slices with the connectedness locus of the two-dimensional uni-Julia sets. We suggested, based on numerical simulations, that a gradual break in connectedness of the uni-J sets occurs in the proximity of the boundary of the equi-M set; a more precise, qualitative description of this transitions requires an analytic approach that is the focus of our future work.

While, as illustrated by the examples considered in this paper, analytic work is quite possible and seems promising in the case of small networks, it is likely that obtaining any useful results for higher-dimensional systems will require different, or additional techniques. We presented two possible approaches, one based on classification, and one based on statistics. We found that network structures which have identical asymptotic dynamics will continue to do so under changes of the quadratic parameter $c$. This is interesting, since it suggests that some information on the long-term outcome in a dynamic network is wired into the architecture, rather than in the node-wise dynamics. We also proposed an alternative, ``average'' view of asymptotic dynamics, counting the structures which produce a certain behavior versus other behaviors. This is a continuation of similar work the authors have carried in continuous-time systems~\cite{radulescu2015nonlinear}. However, similar concepts in discrete networks of quadratic maps are a lot easier to investigate and present well-posed, feasible mathematical problems; the same questions can easily become intractable when using more complicated node dynamics. This speaks in support of using such simple models to begin understanding the behavior of more complex systems, in which direct results are otherwise unreachable. In addition, when using a simple model of quadratic networks, one can put results in the perspective of the long-standing work with single-node iterations, and better understand the mechanisms of transition between a simple system with one operating unit and a complicated dynamic ensemble.

\subsection{Future work}

In this paper, we considered a statistical outlook on classifying asymptotic dynamics in networks with a prescribed architectural property. Another approach to resolving network complexity in a computationally practical way is to reduce the dimensionality of the graph while preserving the dynamics, by collapsing specific sets of nodes to single nodes. For example, as suggested in our prior research, in a graph with communities, rich clubs or strong components (within which the nodes are more tightly connected), it is possible that the dynamics is more robust to changes of structure \emph{within} these modules, and more vulnerable to changes in the coupling \emph{between} the modules. Then, we will investigate the possibility to classify the ensemble dynamics based on simplified representations of the underlying graph, obtained by identifying the robust formations to simple nodes. This can reduce the classification problem to a working framework of much simpler graphs (e.g. trees, cycles), and would also offer a plausible explanation to the preference of natural systems for such hierarchic structures.

One direction in our future work is aimed at investigating a somewhat different temporal coupling scheme for networks, built on principles of random iteration (reminiscent of Markov chains). From each node $j$, there is a probability $p_{jk}$ for the information to travel along the outgoing edge $E_{jk}$ to the adjacent node $k$, so that $z_j$ will be iterated according to the map $z_k(t) \to z_k(t+1)$ attached to that respective node.  This defines a random $n$-dimensional iteration on $(z_1,z_2,...z_n)$. The probabilities $p_{jk}$ are nonzero only when there is an oriented edge connecting $z_j$ and $z_k$. Additionally, the probabilities out of each node (including self cycles) have to add up to one: $\sum_{k=1}^n p_{jk} = 1$. Comerford~\cite{comerford2006hyperbolic} and Sumi~\cite{sumi2001dynamics} have made, for the past ten years, major contributions to the field of random iterations in the one-dimensional case, proving convergence of the Julia sets under random iterations  of hyperbolic polynomial sequences, and describing a phenomenon of cooperation between generating maps as  a factor decreasing the chaos in the overall system~\cite{sumi2011random}. The extension of any of these concepts and results to dynamic networks would be not only mathematically significant, but also of potentially crucial interest to studying networks in the life sciences which may be governed precisely by these rules.

Finally, an extension with potentially high relevance to computational neuroscience would be introducing time and state-dependent edge weights. One of the most fundamental rules in neurobilogy, quantifying the plasticity of brain connections that underlies processes like learning and memory formation, is Hebb's rule. In its most general form, the rule states that the system strengthens connections between neuron/nodes which have correlated (hence potentially causal) activity. One of the simplest historical implementations of Hebb's rule has been to adjust the weight of the each edge by a ``learning'' term proportional to the product of the states of the adjacent nodes, at each iteration step. Then the dynamics of the system of network edge weights becomes as significant as the dynamics of the nodes themselves, with which they are coupled. The weights converge to an attracting state when the network has learned a certain configuration.

\bibliographystyle{plain}
\bibliography{references_Simone}

\clearpage
\renewcommand{\arraystretch}{0.85}


\section*{Appendix A}

\begin{lemma}
The point $c=-3/4$ is not in the equi-M set for the network given by: $z_1 \to z_1^2+c$, $z_2 \to (az_1+z_2)^2+c$, $z_3 \to (z_1+z_2+bz_3)^2+c$, with connectivity weights $a = -1$, $b=-1$.
\end{lemma}

\noindent {\bf Proof.} It is easy to see that the interval $[-3/4,0]$ is invariant under the iteration of the function $z_1 \to z_1^2-3/4$ (since the minimum value of the function $f(z) = z^2-3/4$ for $z \in [-3/4,0]$ is $f(0)=-3/4>-1$, and the maximum value is $f(-3/4)=-3/16>-3/4$. Since $z_1(0)=0 \in [-3/4,0]$, it follows by induction that $z_1 \in [-3/4,0]$ for all iterates, hence the first node is bounded. 

We will show, using induction, that $-3/4 \leq z_2 \leq 0$ also for all iterates. This is satisfied for $z_2(0) =0$. Suppose that $z_2(t) \in [-3/4,0]$ for some $t \geq 0$; we will show that $z_2(t+1)$ is also within this interval. We know that $-3/4 \leq z_1(t) \leq 0$, and $-3/4 \leq z_2(t) \leq 0$, hence $-3/4 \leq -z_1(t)+z_2(t) \leq 3/4$, and $0 \leq (-z_1(t)+z_2(t))^2 \leq 9/16$. Then $z_2(t+1) \in [-3/4, -3/16] \subset [-3/4,0]$, which concludes the induction and shows that the critical $z_2$ is bounded. Moreover, since $z_1,z_2 \in [-3/4,0]$, it follows that $ -1 \leq z_1+z_2+1/2  \leq 1/2$, hence $\lvert z_1+z_2+1/2 \rvert \leq 1$, which we will use below.

It is easy to calculate that the orbit of $z_3$ grows relatively fast for the first portion of the iteration, so that $z_3(8) > 5$. We will use this to show that, in fact, the orbit of the third node escapes to infinity. First notice that, for all iterates (in particular for $t \geq 8$), we have: 
$$\lvert z_3(t+1) \rvert = \lvert (z_1 + z_2-z_3)^2-3/4 \rvert \geq \lvert (z_1+z_2-z_3) \rvert^2 - 3/4$$

\noindent For simplicity, we left out the index for the current iterate (e.g., $z_1$ above represents $z_1(t)$). This further implies that:
$$\sqrt{\lvert z_3(t+1) \rvert +3/4} \geq \lvert(z_1+z_2+1/2)+(-z_3-1/2) \rvert \geq \lvert -z_3-1/2 \rvert - \lvert z_1+z_2+1/2 \rvert$$

\noindent Since $\lvert z_1+z_2+1/2 \rvert \leq 1$, we further have that
$$\sqrt{\lvert z_3(t+1) \rvert +3/4} \geq  \lvert z_3 \rvert -1/2 -1 \geq \lvert z_3 \rvert -3/2$$

\noindent Since $z_3 > 5$, we can square both sides:
$$\lvert z_3(t+1) \rvert \geq (\lvert z_3 \rvert -3/2)^2 -3/4$$

\noindent We want to show that $(\lvert z_3 \rvert -3/2)^2 -3/4 \geq 2 \lvert z_3 \rvert$. Consider the quadratic function 
$$f(\xi) = (\xi-3/2)^2-3/4 -2\xi = \xi^2 -5\xi +3/2,$$ 

\noindent with roots $0<\xi_1 < \xi_2 < 5$. Since $\lvert z_3 \rvert > 5$, it follows that $f(\lvert z_3 \rvert) >0$, hence
$$(1/3 \lvert z_3 \rvert -3/2)^2 -3/4 - 2 \lvert z_3 \rvert > 0$$

\noindent It follows that, as needed
$$\lvert z_3(t+1) \rvert \geq 2\lvert z_3(t) \rvert \text{ for } t \geq 8$$

\noindent In conclusion, the node-wise orbit $z_3$ escapes to infinity.

\hfill $\Box$


\clearpage
\section*{Appendix B}

\begin{figure}[h!]
\begin{center}
\includegraphics[width=\textwidth]{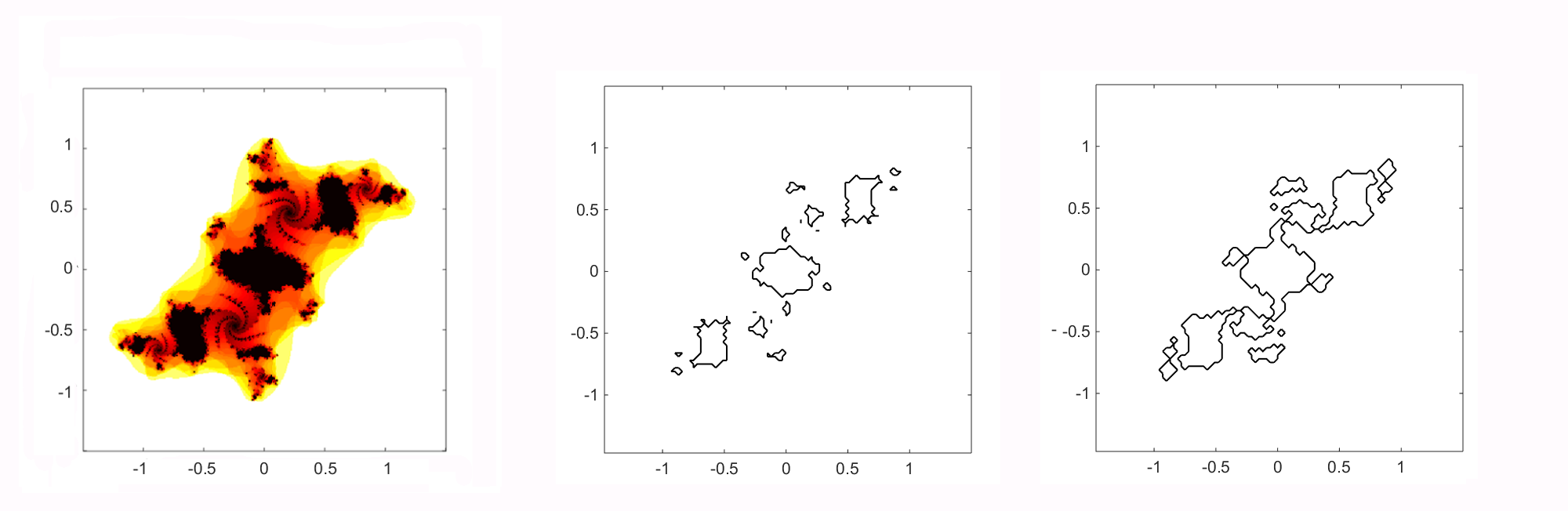}
\end{center}
\caption{\small \emph{{\bf Comparison between detection of connected components of uni-J sets}, using the standard algorithm from the Matlab image processing toolbox versus an improved version including an initial blowup of the Julia set by a one pixel margin. {\bf A.} High resolution ($400 \times 400$ pixels) Uni-J set for the self-drive network $a=-2/3$, $b=-1/3$, corresponding to $c=-0.06-0.68i$; {\bf B.} Count of connected components in low resolution ($100 \times 100$ pixels)  using the standard algorithm found 29 components; {\bf C.} Count of connected components in low resolution ($100 \times 100$ pixels)  using the improved algorithm found 3 components.}}
\label{connect_compare}
\end{figure}

\begin{figure}[h!]
\begin{center}
\includegraphics[width=0.75\textwidth]{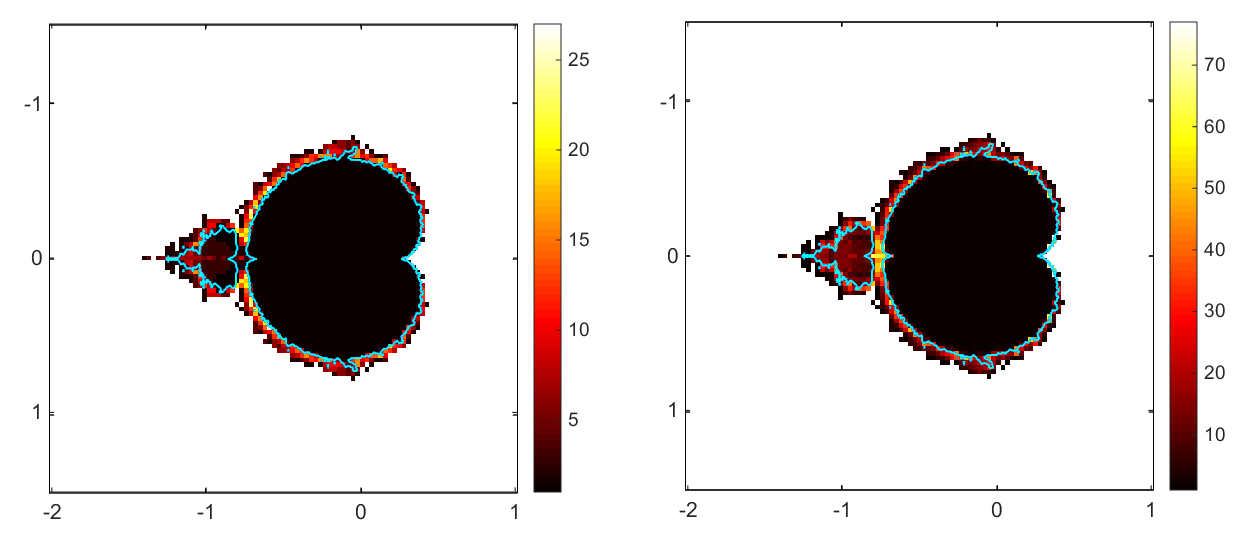}
\end{center}
\caption{\small \emph{{\bf Comparison between the uni-J set connectedness locus} computed using a direct estimate of the number of connected components, versus using the blowup technique. Both panels represent the square $[-2,1] \times [-1.5,1.5]$ in the equi-parameter plane. The blue curve represents the boundary of the equi-M set, computed with 50 iterations. The colors correspond to the number of connected components for the respective uni-J set, computed directly (left) versus using a 1.5 pixel border for the Julia set (right). The panels are almost identical in the black (connected) and white (totally disconnected) regions, while the scale/ number of connected components are very different in the transitional colored region (as shown by the ranges on the color bars).}}
\label{connect_compare}
\end{figure}

\clearpage

\section*{Appendix C}

\begin{figure}[h!]
\begin{minipage}{0.62\textwidth}
\label{table_1_3}
\begin{footnotesize}
\begin{tabular}{|c|c|c|}
\hline
& &\\
$\;\left[ \begin{array}{ccc} 1 & 1 & 1\\ 1 & 1 & 0\\ 1 & 1 & 0 \end{array} \right] \; ({\cal A}_i)$ &
$\;\left[ \begin{array}{ccc} 1 & 1 & 0\\ 1 & 1 & 1\\ 1 & 1 & 0 \end{array} \right] \; ({\cal A}_i)$ &
$\;\left[ \begin{array}{ccc} 1 & 1 & 0\\ 1 & 1 & 0\\ 1 & 1 & 1 \end{array} \right] \; ({\cal B}_{iv})$\\
& &\\
\hline
& &\\
$\;\left[ \begin{array}{ccc} 1 & 1 & 1\\ 1 & 1 & 1\\ 1 & 0 & 0 \end{array} \right] \; ({\cal A}_{ii})$ &
$\;\left[ \begin{array}{ccc} 1 & 1 & 1\\ 1 & 1 & 0\\ 1 & 0 & 1 \end{array} \right] \; ({\cal C}_i)$ &
$\;\left[ \begin{array}{ccc} 1 & 1 & 0\\ 1 & 1 & 1\\ 1 & 0 & 1 \end{array} \right] \; ({\cal D}_v)$ \\
& &\\
\hline
& &\\
$\;\left[ \begin{array}{ccc} 1 & 1 & 1\\ 1 & 0 & 1\\ 1 & 1 & 0 \end{array} \right] \; ({\cal E}_i)$ &
$\;\left[ \begin{array}{ccc} 1 & 1 & 1\\ 1 & 0 & 0\\ 1 & 1 & 1 \end{array} \right] \; ({\cal A}_{ii})$ &
$\;\left[ \begin{array}{ccc} 1 & 1 & 0\\ 1 & 0 & 1\\ 1 & 1 & 1 \end{array} \right] \; ({\cal F}_{iii})$ \\
& &\\
\hline
& &\\
$\;\left[ \begin{array}{ccc} 1 & 1 & 1\\ 1 & 0 & 1\\ 1 & 0 & 1 \end{array} \right] \; ({\cal A}_i)$ &
$\;\left[ \begin{array}{ccc} 1 & 0 & 1\\ 1 & 1 & 1\\ 1 & 1 & 0 \end{array} \right] \; ({\cal F}_{iii})$ &
$\;\left[ \begin{array}{ccc} 1 & 0 & 1\\ 1 & 1 & 0\\ 1 & 1 & 1 \end{array} \right] \; ({\cal D}_v)$ \\
& &\\
\hline
& &\\
$\;\left[ \begin{array}{ccc} 1 & 0 & 0\\ 1 & 1 & 1\\ 1 & 1 & 1 \end{array} \right] \; ({\cal B}_{vi})$ &
$\;\left[ \begin{array}{ccc} 1 & 0 & 1\\ 1 & 1 & 1\\ 1 & 0 & 1 \end{array} \right] \; ({\cal B}_{iv})$ &
$\;\left[ \begin{array}{ccc} 1 & 0 & 1\\ 1 & 0 & 1\\ 1 & 1 & 1 \end{array} \right] \; ({\cal A}_{i})$ \\
& &\\
\hline
& &\\
$\;\left[ \begin{array}{ccc} 1 & 1 & 1\\ 1 & 1 & 1\\ 0 & 1 & 0 \end{array} \right] \; ({\cal A}_{ii})$ &
$\;\left[ \begin{array}{ccc} 1 & 1 & 1\\ 1 & 1 & 0\\ 0 & 1 & 1 \end{array} \right] \; ({\cal D}_{v})$ &
$\;\left[ \begin{array}{ccc} 1 & 1 & 0\\ 1 & 1 & 1\\ 0 & 1 & 1 \end{array} \right] \; ({\cal C}_i)$ \\
& &\\
\hline
& &\\
$\;\left[ \begin{array}{ccc} 1 & 1 & 1\\ 1 & 1 & 1\\ 0 & 0 & 1 \end{array} \right] \; ({\cal B}_{vi})$ &
$\;\left[ \begin{array}{ccc} 1 & 1 & 1\\ 1 & 0 & 1\\ 0 & 1 & 1 \end{array} \right] \; ({\cal F}_{iii})$ &
$\;\left[ \begin{array}{ccc} 1 & 0 & 1\\ 1 & 1 & 1\\ 0 & 1 & 1 \end{array} \right] \; ({\cal D}_{v})$ \\
& &\\
\hline
& &\\
$\;\left[ \begin{array}{ccc} 1 & 1 & 1\\ 0 & 1 & 1\\ 1 & 1 & 0 \end{array} \right] \; ({\cal F}_{iii})$ &
$\;\left[ \begin{array}{ccc} 1 & 1 & 1\\ 0 & 1 & 0\\ 1 & 1 & 1 \end{array} \right] \; ({\cal B}_{vi})$ &
$\;\left[ \begin{array}{ccc} 1 & 1 & 0\\ 0 & 1 & 1\\ 1 & 1 & 1 \end{array} \right] \; ({\cal D}_v)$ \\
& &\\
\hline
& &\\
$\;\left[ \begin{array}{ccc} 1 & 1 & 1\\ 0 & 1 & 1\\ 1 & 0 & 1 \end{array} \right] \; ({\cal D}_{v})$ &
$\;\left[ \begin{array}{ccc} 1 & 1 & 1\\ 0 & 0 & 1\\ 1 & 1 & 1 \end{array} \right] \; ({\cal A}_{ii})$ &
$\;\left[ \begin{array}{ccc} 1 & 0 & 1\\ 0 & 1 & 1\\ 1 & 1 & 1 \end{array} \right] \; ({\cal C}_i)$ \\
& &\\
\hline
& &\\
$\;\left[ \begin{array}{ccc} 1 & 1 & 1\\ 0 & 1 & 1\\ 0 & 1 & 1 \end{array} \right] \; ({\cal B}_{iv})$ &
$\;\left[ \begin{array}{ccc} 0 & 1 & 1\\ 1 & 1 & 1\\ 1 & 1 & 0 \end{array} \right] \; ({\cal E}_i)$ &
$\;\left[ \begin{array}{ccc} 0 & 1 & 1\\ 1 & 1 & 0\\ 1 & 1 & 1 \end{array} \right] \; ({\cal F}_{iii})$ \\
& &\\
\hline
& &\\
$\;\left[ \begin{array}{ccc} 0 & 1 & 0\\ 1 & 1 & 1\\ 1 & 1 & 1 \end{array} \right] \; ({\cal A}_{ii})$ &
$\;\left[ \begin{array}{ccc} 0 & 1 & 1\\ 1 & 1 & 1\\ 1 & 0 & 1 \end{array} \right] \; ({\cal F}_{iii})$ &
$\;\left[ \begin{array}{ccc} 0 & 1 & 1\\ 1 & 0 & 1\\ 1 & 1 & 1 \end{array} \right] \; ({\cal A}_i)$ \\
& &\\
\hline
& &\\
$\;\left[ \begin{array}{ccc} 0 & 0 & 1\\ 1 & 1 & 1\\ 1 & 1 & 1 \end{array} \right] \; ({\cal A}_{ii})$ &
$\;\left[ \begin{array}{ccc} 0 & 1 & 1\\ 1 & 1 & 1\\ 0 & 1 & 1 \end{array} \right] \; ({\cal A}_i)$ &
$\;\left[ \begin{array}{ccc} 0 & 1 & 1\\ 0 & 1 & 1\\ 1 & 1 & 1 \end{array} \right] \; ({\cal A}_i)$ \\
& &\\
\hline
\end{tabular}
\end{footnotesize}
\end{minipage}
\quad
\begin{minipage}{0.395\textwidth}
\includegraphics[width=\textwidth]{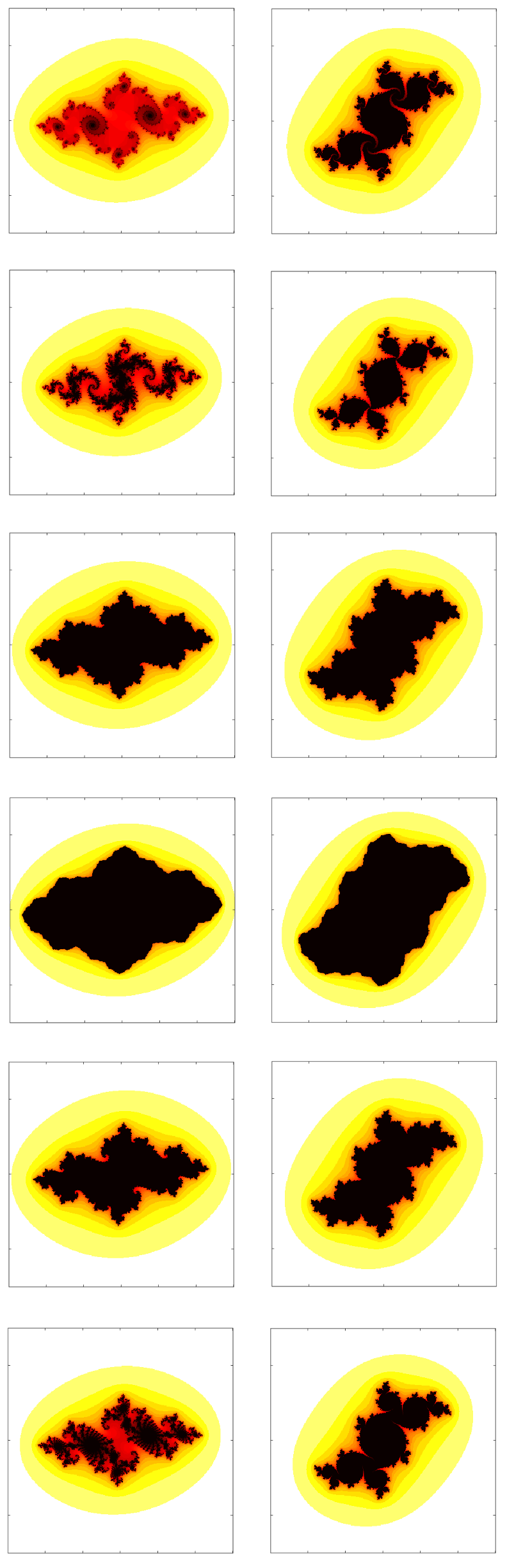}
\end{minipage}
\caption{\label{table_3_3} \emph{\small{\bf Spectral classes versus asymptotic classes} for all networks with $N=3$ nodes and $j=7$ edges. Spectral classes are designated by letters ${\cal A}-{\cal F}$; the asymptotic classes, designate by indices $i$-$vi$, are illustrated on the right for two distinct values of the equi-parameter: $c=-1.15+0.26i$ (left column) and $c=-0.13+i$ (right column). The edge weights were fixed to $g=1/3$. The figure panels show, top to bottom, all asymptotic classes $i$-$vi$, and were created based on 100 iterations, in $400 \times 400$  resolution.}}
\label{classes_simple}
\end{figure}

\end{document}